\documentclass[11pt]{amsart}	

\usepackage{graphicx}
\usepackage{subcaption}
\usepackage{amsmath}
\usepackage{amssymb}
\usepackage{amsthm}
\usepackage{ifthen} 
\usepackage{color}

\usepackage{mathrsfs}

\usepackage[usenames,dvipsnames]{pstricks}
\usepackage{pstricks-add}
\usepackage{epsfig}
\usepackage{pst-grad} 
\usepackage{pst-plot} 
\usepackage[space]{grffile} 
\usepackage{etoolbox} 
\makeatletter 
\patchcmd\Gread@eps{\@inputcheck#1 }{\@inputcheck"#1"\relax}{}{}
\makeatother

\usepackage{dsfont}

\def\polk#1{\setbox0=\hbox{#1}{\ooalign{\hidewidth\lower1.5ex\hbox{`}\hidewidth\crcr\unhbox0}}}

\newcommand{\llip}{\operatorname{Log-Lip}}

\newcommand\restr[2]{{
  \left.\kern-\nulldelimiterspace 
  #1 
  \vphantom{\big|} 
  \right|_{#2} 
  }}

\newtheorem{Theorem}{Theorem}[section]
\newtheorem{Definition}{Definition}[section]
\newtheorem{Lemma}{Lemma}[section]
\newtheorem{Corollary}{Corollary}[section]
\newtheorem{Proposition}{Proposition}[section]
\newtheorem{Condition}{Condition}[section]

\theoremstyle{definition}

\theoremstyle{remark}
\newtheorem{Remark}{Remark}

\allowdisplaybreaks

\newcommand{\intav}[1]{\mathchoice {\mathop{\vrule width 6pt height 3 pt depth -2.5pt\kern -8pt \intop}\nolimits_{\kern -6pt#1}} {\mathop{\vrule width5pt height 3 pt depth -2.6pt \kern -6pt \intop}\nolimits_{#1}}{\mathop{\vrule width 5pt height 3 pt depth -2.6pt \kern -6pt\intop}\nolimits_{#1}} {\mathop{\vrule width 5pt height 3 pt depth-2.6pt \kern -6pt \intop}\nolimits_{#1}}}

\usepackage{setspace}
\numberwithin{equation}{section}

\title[Transmission problems]{Transmission problems: regularity theory, interfaces and beyond}

\author[V. Bianca]{Vincenzo Bianca}
\address{University of Coimbra, CMUC, Department of Mathematics, 3001-501 Coimbra, Portugal}{}
\email{vincenzo@mat.uc.pt}

\author[E. A. Pimentel]{Edgard A. Pimentel}
\address{University of Coimbra, CMUC, Department of Mathematics, 3001-501 Coimbra, Portugal and Pontifical Catholic University of Rio de Janeiro -- PUC-Rio, 22451-900 G\'avea, Rio de Janeiro-RJ, Brazil}{}
\email{edgard.pimentel@mat.uc.pt}

\author[J.M.~Urbano]{Jos\'{e} Miguel Urbano}
\address{King Abdullah University of Science and Technology (KAUST), Computer, Electrical and Mathematical Sciences and Engineering Division (CEMSE), Thuwal 23955-6900, Saudi Arabia and University of Coimbra, CMUC, Department of Mathematics, 3001-501 Coimbra, Portugal}{} 
\email{miguel.urbano@kaust.edu.sa} 

\date{\today}

\begin{document}

\begin{abstract}
Modelling diffusion processes in heterogeneous media requires addressing inherent discontinuities across interfaces, where specific conditions are to be met. These challenges fall under the purview of Mathematical Analysis as \emph{transmission problems}. We present a pa\-no\-ra\-ma of the theory of transmission problems, encompassing the seminal contributions from the 1950s and subsequent developments. Then we delve into the discussion of regularity issues, including recent advances matching the minimal regularity requirements of interfaces and the optimal regularity of the solutions. A discussion on free transmission problems closes the survey.
\end{abstract}



\maketitle

\section{Introduction}\label{sec_philiproth}


The study of the physical properties of heterogeneous media plays a fundamental role across disciplines, both in the basic sciences and in their applications. An underlying phenomenon governing an important set of those properties is \emph{diffusion}. Given the heterogeneity of the media, it is reasonable to consider the diffusion process changes across the different components of the medium. Since the second half of the last century, methods in the realm of partial differential equations (PDE, for short) and the calculus of variations have been used to examine these changes. Their various mathematical formulations are known in the literature as \emph{transmission problems}.

The general setting is the following. Let $\Omega\subseteq\mathbb{R}^d$ be a bounded open set; for $k\in\mathbb{N}$ fixed, consider $\Omega_1,\dots,\Omega_k\Subset\Omega$, pairwise disjoint. Suppose $u:\Omega\to\mathbb{R}$ solves within every $\Omega_i$, $i=1,\dots,k$, some prescribed PDE. A natural question concerns the interactions of the $k$ restrictions
\[
    u_i:=\restr{u}{\Omega_i},\hspace{.3in}\mbox{for}\hspace{.3in}i=1,\ldots,k,
\]
with each other. 

At this point, an important distinction is in place. In fact, one must distinguish the case of an interface that is fixed \textit{a priori} from the case of a solution-dependent interface. The latter is referred to as \emph{free interface} or \emph{free boundary}. From a chronological viewpoint, fixed transmission problems precede the free interface variant. 

Transmission problems first appeared in Mathematical Analysis in the work of Mauro Picone, circa 1950 \cite{Picone1954}. The author addresses a fixed transmission problem arising in elasticity theory in his seminal contribution. The main result in \cite{Picone1954} is the uniqueness of solutions for the transmission problem. The author also indicates a strategy aiming at establishing existence, which however is not explored at length in the paper.

The work of Mauro Picone is consequential, as it attracted immediate attention from the mathematical community to this class of problems. We highlight the work of Jacques-Louis Lions \cite{Lions1956}, where the first rigorous theory of the existence of solutions is put forward. Among the subsequent studies inspired by Picone's ideas, we mention \cite{Lions-Schwartz1955, Stampacchia1956, Campanato1957, Campanato1959, Campanato1959a, Iliin-Shismarev1961,Oleinik1961, Borsuk1968}. We refer the reader to the neatly written monograph \cite{Borsuk2010}, where a detailed account of the theory of fixed transmission problems and a comprehensive list of references are available. 

Once the questions on the existence and uniqueness of solutions have been well-understood and documented, the attention of the community shifted to the regularity theory. A pioneering contribution in this regard is in the work of Yanyan Li and Michael Vogelius \cite{Li-Vogelius2000}. In that paper, the authors consider a domain $\Omega\subseteq\mathbb{R}^d$, split into $k\in\mathbb{N}$ pairwise disjoint subdomains $\Omega_1,\dots,\Omega_k\Subset\Omega$. The mathematical model is driven by the PDE
\begin{equation*}
    \frac{\partial}{\partial x_i}\left(a(x)\frac{\partial u}{\partial x_j}\right)=f\hspace{.2in}\mbox{in}\hspace{.2in}\Omega,
\end{equation*}
where $f:\Omega\to\mathbb{R}$ satisfies usual continuity and integrability assumptions, and 
\[
a(x)=
\begin{cases}
a_i(x)&\hspace{.1in}\mbox{for}\hspace{.1in}\Omega_i,\,i=1,\dots,k\\
a_{k+1}(x)&\hspace{.1in}\mbox{\textnormal{for}}\hspace{.1in}\Omega\setminus\cup_{i=1}^k\Omega_i.
\end{cases}
\]
The main contribution of \cite{Li-Vogelius2000} is the local H\"older-continuity for the gradient of the solutions. The vectorial counterpart of the results in \cite{Li-Vogelius2000} appeared in the work of Yanyan Li and Louis Nirenberg \cite{Li-Nirenberg2003}, where the authors also establish estimates for higher-order derivatives of the solutions. From the perspective of applications, this model accounts for different inclusions within a homogeneous medium, the primary example being fibre-reinforced materials.

The regularity analysis in the context of conductivity and insulation is the subject of \cite{Bao-Li-Yin1, Bao-Li-Yin2}. Here, the authors consider a bounded domain $\Omega\subseteq\mathbb{R}^d$, and two disjoint subdomains $\Omega_1,\Omega_2\Subset\Omega$ that are $\varepsilon$-apart, that is $\textnormal{dist}(\Omega_1,\Omega_2)=\varepsilon$. Within each subdomain, we are in the presence of an equation in divergence form. The analysis of the problem leads to estimates which are dependent on the parameter $\varepsilon>0$. Indeed, they ensure a deterioration on the bounds for the gradient, which blows up as $\varepsilon\to0^+$.

In the former literature, the program has focused on the existence and regularity of the solutions under regular enough fixed interfaces. Recently, the fundamental question of examining the \emph{minimal} geometric requirements on the fixed interface to recover regularity estimates on the solutions of the fixed transmission problem has been addressed. This is the analysis in the work of Luis Caffarelli, Mar\'ia Soria-Carro and Pablo Stinga \cite{CSCS2021}, where the authors consider a bounded domain $\Omega\subseteq\mathbb{R}^d$, and fix $\Omega_1\Subset\Omega$, defining $\Omega_2:=\Omega\setminus\Omega_1$. The key assumption in \cite{CSCS2021} concerns the geometry of $\Gamma:=\partial\Omega_1$, which the authors assume to be of class $C^{1,\alpha}$. 

The problem is to look for a function $u\in W^{1,2}_{\rm loc}(\Omega)$ whose restrictions $u_i:=u\big|_{\overline{\Omega_i}}$ satisfy
\begin{equation}\label{eq_oxford}
\begin{cases}
\Delta u_i=0&\hspace{.1in}\mbox{\textnormal{in}}\hspace{.1in}\Omega_i,\hspace{.1in}\mbox{for}\hspace{.1in}i=1,2,\\
u_2=0&\hspace{.1in}\mbox{\textnormal{in}}\hspace{.1in}\partial\Omega,\\
(u_1)_\nu-(u_2)_\nu=g&\hspace{.1in}\mbox{\textnormal{on}}\hspace{.1in}\Gamma.
\end{cases}
\end{equation}

By resorting to well-known properties of harmonic functions, the authors prove the existence and uniqueness of the solutions, together with a regularity result in $C^{0,{\rm Log-Lip}}$-spaces \emph{across} the fixed interface $\Gamma$. 

However, the main contribution is \cite{CSCS2021} concerns the regularity of solutions \emph{up to the interface}. The authors develop a new stability result relating general $C^{1,\alpha}$-regular interfaces with flat ones. Their approach relies on a clever combination of the mean value formula and the maximum principle for harmonic functions. As a result, they show that $u_i\in C_{\rm loc}^{1,\alpha}\left(\overline{\Omega_i}\right)$ for $i=1,2$. 

In \cite{SoriaCarro-Stinga} the authors extend the analysis put forward in \cite{CSCS2021} to the context of fully nonlinear equations. In that paper, a complete theory of viscosity solutions to fully nonlinear fixed transmission problems is developed, including existence and uniqueness, comparison principles, and regularity estimates. As in \cite{CSCS2021}, the authors establish the regularity of the solutions up to the fixed interface; we notice the regularity regime for the solutions matches the geometric properties of the fixed interface. Regarding the $C^{1,\alpha}$-regularity of solutions up to the interface, we refer the reader to \cite{DongP} for a variant of the proof put forward in \cite{CSCS2021}. The argument in \cite{DongP} does not rely either on the mean value property or on the maximum principle for harmonic functions. The paper also includes extensions to the context of $C^{1,{\rm Dini}}$-fixed interfaces.

This introduction has covered so far some of the developments in the context of \emph{fixed transmission problem}. That is, in the presence of interfaces that are determined \textit{a priori}. An important variant of the topic regards the so-called \emph{free transmission problems}, where the interface depends on the solution of the problem. The first attempt to examine a free transmission problem appeared in \cite{Amaral-Teixeira2015}, where the authors propose a variational model. More precisely, they consider maps $A_\pm:\Omega\to S(d)$, where $S(d)$ stands for the space of symmetric matrices of order $d$, and study the functional 
\begin{equation}\label{eq_roxo}
    I(v):=\int_{\Omega}\frac{1}{2}\left\langle A(v,x)Dv,Dv\right\rangle+\Lambda(v,x)+fv\,{\rm d}x,
\end{equation}
where the matrix $A(v,x)$ and the function $\Lambda(v,x)$ are defined as
\[
    A(v,x):= A_+(x)\chi_{\{v>0\}}+A_-(x)\chi_{\{v<0\}}
\]
and
\[
    \Lambda(v,x):=\lambda_+(x)\chi_{\{v>0\}}+\lambda_-(x)\chi_{\{v<0\}},
\]
and $f\in L^\infty(\Omega)$. The program pursued in \cite{Amaral-Teixeira2015} includes the existence of minimizers for \eqref{eq_roxo}, together with regularity results in the spaces of H\"older and Lipschitz-continuous functions.

In the context of free transmission problems at the intersection of PDE and the calculus of variations, we mention the recent work of Maria Colombo, Sungham Kim and Henrik Shahgholian \cite{Colombo-Kim-Shahgholian2023}. The authors consider the functional
\[
	J(v,\Omega):=\int_\Omega\left(|Dv^+|^p+|Dv^-|^q\right){\rm d}x,
\]
where $1<p,q<\infty$, with $p\neq q$. Their findings include properties of minimizers as well as information on the associated free boundary. They prove the existence of minimizers, and their local H\"older continuity, with estimates. Concerning the analysis of the free boundary, the authors derive a free boundary condition and prove the free boundary is of class $C^{1,\alpha}$, almost everywhere, with respect to the measure $\Delta_pu^+$. The paper also verifies the support of $\Delta_pu^+$ is of $\sigma$-finite $(n-1)$-dimensional Hausdorff measure. 

In the realm of fully nonlinear operators, free transmission problems have been studied from two different perspectives, namely the existence of generalized solutions (in the viscosity and the strong senses), and their regularity properties (see \cite{Pimentel-Swiech2022,Huaroto-Pimentel-Rampasso-Swiech}). In the uniformly elliptic problem, the model under consideration has the form
\begin{equation}\label{eq_scorpion}
    F_1(D^2u)\chi_{\{u>0\}}+F_2(D^2u)\chi_{\{u<0\}}=f\hspace{.15in}\mbox{in}\hspace{.15in}\Omega\cap\left(\left\lbrace u>0\right\rbrace\cup \left\lbrace u<0\right\rbrace\right),
\end{equation}
where $F_i:S(d)\to \mathbb{R}$ are fully nonlinear uniformly elliptic operators and $f\in L^\infty(\Omega)$ is a given source term.

The existence of solutions for the Dirichlet problem associated with \eqref{eq_scorpion}, both in the $L^p$-viscosity and in the $L^p$-strong sense, was established in \cite{Pimentel-Swiech2022}. The regularity of the solutions is the subject of \cite{Pimentel-Santos2023}. 

A degenerate variant of \eqref{eq_scorpion} has also been examined in the literature. Let $0<p_1<p_2$ be fixed positive numbers and consider
\begin{equation}\label{eq_policeextorsion}
    |Du|^{p_1\chi_{\{u>0\}}+p_2\chi_{\{u<0\}}}F(D^2u)=f\hspace{.15in}\mbox{in}\hspace{.15in}\Omega\cap\left(\left\lbrace u>0\right\rbrace\cup \left\lbrace u<0\right\rbrace\right).
\end{equation}
The model in \eqref{eq_policeextorsion} is inspired by the tradition of fully nonlinear equations degenerating as a power of the gradient; see \cite{Birindelli-Demengel2004,Birindelli-Demengel2006,Birindelli-Demengel2007a,Birindelli-Demengel2007b,Birindelli-Demengel2009,Birindelli-Demengel,Imbert-Silvestre,Araujo-Ricarte-Teixeira}, among others. In \cite{Huaroto-Pimentel-Rampasso-Swiech}, the authors prove the existence of solutions to \eqref{eq_policeextorsion} and examine their regularity properties. In particular, they prove a $C^{1,\alpha}$-regularity estimate. Finally, the modulus of continuity depends explicitly on the Krylov-Safonov theory available for $F=0$ and on the degeneracy rates $0<p_1<p_2$. 

The variant of \eqref{eq_policeextorsion} with nonhomogeneous degeneracies is the subject of \cite{CDF2022a}. In that paper, the author considers a model of the form
\[
    H(x,u,Du)F(D^2u)=f\hspace{.2in}\mbox{in}\hspace{.2in}\Omega,
\]
where $F$ is uniformly elliptic, $f\in C(\Omega)\cap L^\infty(\Omega)$, and 
\[
    H(x,u,Du):=|Du|^{p_1\chi_{\{u>0\}}+p_2\chi_{\{u<0\}}}+a(x)\chi_{\{u>0\}}|Du|^q+b(x)\chi_{\{u<0\}}|Du|^s,
\]
for nonnegative functions $a,b\in C(\Omega)$, and nonnegative exponents $p_1$, $p_2$, $q$, and $s$. The findings in \cite{CDF2022a} cover the existence of solutions and $C^{1,\alpha}$-regularity estimates.

The remainder of this paper is organized as follows. Section \ref{sec_abelprize} discusses a few examples of fixed transmission problems, covering some of the main aspects of the theory. Section \ref{subsec_columbus} describes the seminal work of Mauro Picone \cite{Picone1954} and presents an argument due to Jacques-Louis Lions for the existence of weak solutions to Picone's problem. More recent results on regularity theory are the subject of Section \ref{subsec_tar}. In Section \ref{subsec_tango} we examine the regularity theory of fixed transmission problems in the presence of minimal geometric requirements on the transmission interface. We describe quasilinear degenerate problems in Section \ref{sec_cpr}, whereas their regularity in borderline spaces is the subject of Section \ref{subsec_bmo}. The paper ends with a discussion on free transmission problems in Section \ref{subsec_ussr}.

\section{An overview on fixed transmission problems}\label{sec_abelprize}

\subsection{Laying the foundations: first models and a theory of existence}\label{subsec_columbus}

In this section, we describe the foundations of the theory of transmission problems in detail. In $\mathbb{R}^d$, with $d\ge3$, let $\Omega_1$ and $\Omega_2$ be two disjoint open sets. Suppose that $\Omega_2$ is the complement of a compact set. Suppose that $\partial\Omega_1$ and $\partial\Omega_2$ share an $(d-1)$-dimensional manifold $\Gamma$, continuously differentiable except, at most, at a finite set of points. For convenience, set $S_i:=\partial\Omega_i\setminus\Gamma$, for $i=1,2$, and
\begin{equation}\label{Picone}
    \Omega:=\Omega_1\cup\Omega_2.
\end{equation}
Now, let $h_1,k_1$ (resp. $h_2,k_2$) be the constants of Lamé of $\Omega_1$ (resp. $\Omega_2$). Lamé constants are two material-dependent quantities that arise in strain-stress relationships. These parameters are named after Gabriel Lamé\footnote{Gabriel Lam\'e was a French mathematician, born in Tours in the year 1795. With various contributions to a wide range of topics in Mathematics, he is known for a theory of curvilinear coordinates and the study of ellipse-like curves, currently referred to as \emph{Lam\'e curves}. A Foreign Member of the Royal Swedish Academy of Sciences, elected in 1854, Lam\'e is one of the 72 names inscribed in the Eiffel Tower. He died in Paris, in 1870.}. One looks for two vector fields, $U^1$ and $U^2$, defined in $\Omega_1$, and $\Omega_2$ respectively, with values in $\mathbb{R}^d$, such that
\begin{equation}\label{Picone_problem}
\begin{cases}
h_1\Delta U^1+(h_1+k_1)D({\rm div}\,U^1)+f^1=0&\hspace{.1in}\mbox{\textnormal{in}}\hspace{.1in}\Omega_1\\
h_2\Delta U^2+(h_2+k_2)D({\rm div}\,U^2)+f^2=0&\hspace{.1in}\mbox{\textnormal{in}}\hspace{.1in}\Omega_2,
\end{cases}
\end{equation}
under the interface conditions
\begin{equation}\label{Picone_condition}
\begin{cases}
U^1=U^2&\hspace{.1in}\mbox{\textnormal{on}}\hspace{.1in}\Gamma\\
t_1(U^1)+t_2(U^2)=0&\hspace{.1in}\mbox{\textnormal{on}}\hspace{.1in}\Gamma,
\end{cases}
\end{equation}
where $D({\rm div}\,U^i)$ is the gradient of the divergence of $U^i$ for $i=1,2$.

Here, $f^1$ and $f^2$ are two given vector fields, representing external forces, and $t_1(U^1),t_2(U^2)$ are the pressure fields on $\partial\Omega_1$ and $\partial\Omega_2$, respectively. Moreover, $t_i(U^i)$ has to be equal to the vector field on $S_i$, for $i=1,2$, and $U^2(x)\to0$ as $|x|\to\infty$. In \cite{Picone1954}, Picone proved the uniqueness of the solution to this problem, claiming that solutions exist. Our choice in this survey is to present the general lines of an existence theory for \eqref{Picone_problem}-\eqref{Picone_condition}, due to Jacques-Louis Lions \cite{Lions1956}.

We need some ingredients from the theory of functional spaces to discuss the argument put forward by Jacques-Louis Lions leading to the existence of a solution to \eqref{Picone_problem}-\eqref{Picone_condition}. Let $\Omega\subseteq\mathbb{R}^d$ be bounded and open. For $u\in W^{1,2}(\Omega;\mathbb{R}^d)$, with $k,\ell\in\{1,\dots,d\}$, define
\begin{equation*}
s_{k\ell}(u):=\frac{1}{2}\left(\frac{\partial u_k}{\partial x_{\ell}}-\frac{\partial u_\ell}{\partial x_k}\right)
\end{equation*}
and
\begin{equation*}
S(u):=\sum_{k,\ell=1}^d\|s_{k\ell} (u)\|_{L^2(\Omega)}^2.
\end{equation*}
For $u,v\in L^2(\Omega;\mathbb{R}^d)$, we denote by $(u,v)_{L^2(\Omega)}$ the inner product between $u$ and $v$ in $L^2(\Omega;\mathbb{R}^d)$, i.e.,
\begin{equation*}
(u,v)_{L^2(\Omega)}:=\int_\Omega u\cdot v\,{\rm d}x
\end{equation*}
For $u,v\in W^{1,2}(\Omega;\mathbb{R}^d)$, we then set
\begin{equation*}
S(u,v):=\sum_{k,\ell=1}^d\big(s_{k\ell}(u),s_{k\ell}(v)\big)_{L^2(\Omega)}.
\end{equation*}
Note that $S(u)=S(u,u)$, and also
\begin{equation}\label{point_essentiel}
S(u)\le\sum_{k,\ell=1}^d\left\|\frac{\partial u_k}{\partial x_\ell}\right\|_{L^2(\Omega)}^2.
\end{equation}

The definition of the following space will play a central role in the sequel. We recall the Sobolev exponent $2^*>2$ is given by
\begin{equation*}
    2^*:=\frac{2d}{d-2},
\end{equation*}
whereas its H\"older conjugate is
\[
    (2^*)':=\frac{2d}{d+2}.
\]

\begin{Definition}
Let $\Omega$ be as in \eqref{Picone}. For $i=1,2$, let $T:W^{1,2}(\Omega_i;\mathbb{R}^d)\to L^2(\Gamma;\mathbb{R}^d)$ be the trace operator. We define $V$ to be the space of functions belonging to $L^{2^*}(\Omega;\mathbb{R}^d)\cap W^{1,2}(\Omega;\mathbb{R}^d)$ such that $T(u^1)=T(u^2)$, where $u^1\in W^{1,2}(\Omega_1;\mathbb{R}^d)$ and $u^2\in W^{1,2}(\Omega_2;\mathbb{R}^d)$.
\end{Definition}

One fundamental step in the reasoning put forward in \cite{Lions1956} is to prove the inequality in \eqref{point_essentiel} can be reversed, up to constant factors, in $V$. Thanks to this, $V$ will be a Hilbert space endowed with a natural norm, and the existence of an isomorphism between a suitable subspace of $V$ and $L^\frac{2d}{d+2}(\Omega;\mathbb{R}^d)$ will lead to the existence of a unique solution. Before stating a result on the existence of solutions to \eqref{Picone_problem}-\eqref{Picone_condition}, we detail two definitions concerning open and connected subsets of $\mathbb{R}^d$.  

\begin{Definition}[Sobolev set]\label{def_sobolev}
Let $\Omega\subseteq\mathbb{R}^d$ be connected and open. If $|\Omega|<\infty$, we say that $\Omega$ is a Sobolev set if $W^{1,2}(\Omega)\subset L^{2^*}(\Omega)$. If $|\Omega|=\infty$, we say that $\Omega$ is a Sobolev set if for every $u\in W^{1,2}(\Omega)$ there exists a constant $c=c(u)$ such that $u+c(u)\in L^{2^*}(\Omega)$.
\end{Definition}

\begin{Remark}[Sobolev sets and the geometry of domains]\label{rem_vivaldi}
We notice that if the boundary of $\Omega$ is locally Lipschitz-regular, the inclusion $W^{1,2}(\Omega)\subset L^{2^*}(\Omega)$ follows from standard results in the theory of Sobolev spaces. As a consequence, every Lipschitz-regular domain is a Sobolev set. However, the notion of \emph{Sobolev set} also accommodates less regular subsets of $\mathbb{R}^d$ for which that inclusion is not available.
\end{Remark}

\begin{Definition}[Friedrichs set]\label{def_coffee}
Let $\Omega\subseteq\mathbb{R}^d$ be connected and open. We say that $\Omega$ is a Friedrichs set if there exist $K\subseteq\Omega$, compact, and a constant $C\in(0,1)$ such that, for every $u\in W^{1,2}(\Omega;\mathbb{R}^d)$ satisfying
\begin{equation*}
-\Delta u-D({\rm div}\,u)=0\hspace{.3in}\mbox{in}\hspace{.3in}\Omega,
\end{equation*}
we have
\begin{equation*}
    S(u)\le \frac{1}{1-C}S_K(u),
\end{equation*}
where 
\[
    S_K(u):=\sum_{k,\ell=1}^d\|s_{k\ell}\|_{L^2(K)}^2.
\]
\end{Definition}
We refer the reader to \cite{Friedrichs1,Friedrichs2} for further context on Friedrichs sets. The main result in \cite{Lions1956} is the following theorem.

\begin{Theorem}[Existence of a unique solution]\label{existence_Picone}
Let $\Omega\subseteq\mathbb{R}^d$ be as in \eqref{Picone}. Suppose that $\Omega_{1}$ and $\Omega_{2}$ are Sobolev and Friedrichs sets.
Then there exists a unique solution to the problem of Picone \eqref{Picone_problem}-\eqref{Picone_condition}.
\end{Theorem}

Now, we continue by detailing the general lines of the proof of Theorem \ref{existence_Picone}. We proceed with a lemma related to the space of tempered distributions. First, equip the linear space $L^{2^*}(\Omega;\mathbb{R}^d)\cap W^{1,2}(\Omega;\mathbb{R}^d)$ with the norm
\[
    \left\|\,\cdot\,\right\|_{2^*,2,\Omega}:L^{2^*}(\Omega;\mathbb{R}^d)\cap W^{1,2}(\Omega;\mathbb{R}^d)\to[0,\infty),
\]
given by 
\[
    \left\|u\right\|_{2^*,2,\Omega}:=\left\|u\right\|_{L^{2^*}(\Omega;\,\mathbb{R}^d)}+\left\|Du\right\|_{L^2(\Omega;\,\mathbb{R}^d)}.
\]
We denote with $\mathcal{D}_{2^*,2}(\Omega;\mathbb{R}^d)$ the closure of $\mathcal{D}(\Omega;\mathbb{R}^d)$ in the space $L^{2^*}(\Omega;\mathbb{R}^d)\cap W^{1,2}(\Omega;\mathbb{R}^d)$, with respect to the norm $\left\|\,\cdot\,\right\|_{2^*,2,\Omega}$. 

\begin{Lemma}\label{lem_easyjet}
For every $u\in\mathcal{D}_{2^*,2}(\Omega;\mathbb{R}^d)$ we have
\begin{equation*}
\sum_{k,\ell=1}^d\left\|\frac{\partial u_k}{\partial x_\ell}\right\|_{L^2(\Omega)}^2\le2S(u).
\end{equation*}
\end{Lemma}

The previous result, whose proof follows by a simple integration by parts and by resorting to a density argument, says that the inequality in \eqref{point_essentiel} can be reversed in $\mathcal{D}_{2^*,2}(\Omega;\mathbb{R}^d)$. This is a preliminary step in proving one can reverse \eqref{point_essentiel} in the space $V$.

Define the bilinear form $\left\langle\cdot,\cdot\right\rangle:W^{1,2}(\Omega;\mathbb{R}^d)\to\mathbb{R}$ as
\begin{equation*}
\left\langle u,v\right\rangle:=2hS(u,v)+k({\rm div}\,u,{\rm div}\,v)_{L^2(\Omega)}
\end{equation*}
with $h>0$ and $3k+2h>0$. The operator $L$ defined as
\begin{equation}\label{L}
L:=-h\Delta-(h+k)D{\rm div}
\end{equation}
is such that
\begin{equation*}
\left\langle u,v\right\rangle=\left\langle Lu,v\right\rangle_{(\mathcal{D}_{2^*,2}(\Omega;\mathbb{R}^d)',\mathcal{D}_{2^*,2}(\Omega;\mathbb{R}^d))},
\end{equation*}
where the crochet at the right-hand side represents the duality between $\mathcal{D}_{2^*,2}(\Omega;\mathbb{R}^d)'$ and $\mathcal{D}_{2^*,2}(\Omega;\mathbb{R}^d)$. The operator \eqref{L} is elliptic, \textit{i.e.}, there exists $a>0$ such that
\begin{equation*}
\left\langle u,u\right\rangle\ge aS(u),
\end{equation*}
for every $u\in\mathcal{D}(\Omega;\mathbb{R}^d)$.

Because the operator in \eqref{L} is elliptic, it defines an isomorphism between $\mathcal{D}_{2^*,2}(\Omega;\mathbb{R}^d)'$ and $\mathcal{D}_{2^*,2}(\Omega;\mathbb{R}^d)$. An important ingredient in the proof of Theorem \ref{existence_Picone} concerns extending Lemma \ref{lem_easyjet} to functions in the functional space $V$. This is the content of the next proposition.

\begin{Proposition}\label{prop_avola }
Let $\Omega\subseteq\mathbb{R}^d$ be as in \eqref{Picone}. Suppose that $\Omega_{1}$ and $\Omega_{2}$ are both Sobolev and Friedrichs sets. Then, there exists $C>0$ such that
\begin{equation*}
\sum_{k,\ell=1}^d\left\|\frac{\partial u_k}{\partial x_\ell}\right\|_{L^2(\Omega)}^2\le CS(u),
\end{equation*}
for every $u\in V$.
\end{Proposition}
As an immediate consequence, one has the following corollary.
\begin{Corollary}
Let $\Omega\subseteq\mathbb{R}^d$ be as in \eqref{Picone}. Suppose that $\Omega_{1}$ and $\Omega_{2}$ are both Sobolev and Friedrichs sets. Then $V$ is a Hilbert space endowed with the norm $\sqrt{S(u)}$.
\end{Corollary}

Now that the aspects of functional analysis have been formulated, we rigorously frame Picone's problem \eqref{Picone_problem}-\eqref{Picone_condition} in the context of those elements. We start by defining operators and suitable functional spaces which allow us to produce an existence result.

\subsubsection{An existence result for Picone's problem \eqref{Picone_problem}-\eqref{Picone_condition}}\label{subsubsec_picone} 

Let $\Omega\subseteq\mathbb{R}^d$ be as in \eqref{Picone}. Let $F^i:W^{1,2}(\Omega;\mathbb{R}^d)\times W^{1,2}(\Omega;\mathbb{R}^d)\to\mathbb{R}$, for $i=1,2$, be defined as
\begin{equation*}
F^1(u^1,v^1):=2h_1S(u^1,v^1)+k_1({\rm div}\,u^1,{\rm div}\,v^1)_{L^2(\Omega_{1})},
\end{equation*}
\begin{equation*}
F^2(u^2,v^2):=2h_2S(u^2,v^2)+k_2({\rm div}\,u^2,{\rm div}\,v^2)_{L^2(\Omega_{2})},
\end{equation*}
and
\begin{equation*}
F(u,v):=F^1(u^1,v^1)+F^2(u^2,v^2),
\end{equation*}
where $h_i$ and $k_i$ are the Lamé constants for $\Omega_i$, for $i=1,2$, $u=(u^1,u^2)$ and $v=(v^1,v^2)$.

Let $R_i:L^2(\Gamma)\to L^2(\Gamma)$ be a bounded linear operator for $i=1,2$. Define
\begin{equation*}
G^i(u^i,v^i):=F^i(u^i,v^i)+\big(R_iT(u^i),T(v^i)\big)_{L^2(\Gamma;\mathbb{R}^d)},
\end{equation*}
and
\begin{equation}\label{sesq}
G(u,v):=G^1(u^1,v^1)+G^2(u^2,v^2).
\end{equation}
It is easy to see that \eqref{sesq} is a bilinear continuous operator in $W^{1,2}(\Omega;\mathbb{R}^d)\times W^{1,2}(\Omega;\mathbb{R}^d)$. Moreover, if $R_i$ are positive and small enough, the operator in \eqref{sesq} is elliptic in $V$.

Now that the operators are defined, we can set the functional spaces where to look for the solution $U$ to \eqref{Picone_problem} and verify the transmission conditions \eqref{Picone_condition}.

Define the space $\mathcal{H}$ is the space of functions $u\in W^{1,2}(\Omega;\mathbb{R}^d)$ such that $Lu\in L^{\frac{2d}{d+2}}(\Omega;\mathbb{R}^d)$. Also, define $N$ as the space of functions $u\in V$ such that $Lu\in L^{\frac{2d}{d+2}}(\Omega;\mathbb{R}^d)$, and
\begin{equation*}
\left\langle Lu,v\right\rangle_{(L^{2d/(d+2)}(\Omega;\mathbb{R}^d),L^{2^*}(\Omega;\mathbb{R}^d))}=F(u,v),
\end{equation*}
for every $v\in V$. We endow $N$ with the norm given by
\begin{equation*}
    \|u\|_N:=\|u\|_V+\|Lu\|_{L^{2d/(d+2)}(\Omega;\mathbb{R}^d)}.
\end{equation*}

In what follows, we state a pivotal result from \cite{Lions1956}. 

\begin{Proposition}\label{prop_cool}
Let $\Omega\subseteq\mathbb{R}^d$ be as in \eqref{Picone}. Suppose that $\Omega_{1}$ and $\Omega_{2}$ are both Sobolev and Friedrichs sets. Suppose further that $F$ is elliptic in $V$. Then $L$ is an isomorphism between $N$ and $L^{\frac{2d}{d+2}}(\Omega;\mathbb{R}^d)$.
\end{Proposition}

Briefly, the proof of Proposition \ref{prop_cool} relies on two facts. First, we have that $(V,\sqrt{S(\,\cdot\,)})$ is a Hilbert space. Secondly, one notices that, for every $u\in V$, one can define a linear continuous operator $\Tilde{F}:V\to V$ such that
    \begin{equation*}
        F^2(u,v)=F^1\left(\Tilde{F}(u),v\right),\quad\Tilde{F}(u)\in V,
    \end{equation*}
and eventually $\Tilde{F}$ is a symmetric operator.
    
The proof of Theorem \ref{existence_Picone} is now an easy consequence of the previous proposition; one argues as follows. Start by finding $U\in\mathcal{H}$ such that, for a given $H\in L^{\frac{2d}{d+2}}(\Omega;\mathbb{R}^d)$ and $h\in\mathcal{H}$, we have
\begin{equation*}
LU=H
\end{equation*}
and
\begin{equation}\label{limits}
U-h\in N.
\end{equation}

Now, we give a precise formulation of \eqref{limits} that will lead to the transmission condition stated at the beginning. Integrating by parts, we get
\begin{equation*}
\left\langle L_1u^1,v^1\right\rangle=F(u^1,v^1)-\int_{\partial\Omega_{1}}t_1(u^1)\cdot v^1\,{\rm d}\mathcal{H}^{d-1}.
\end{equation*}
Similarly, one has the formula for $\Omega_{2}$. Now, $u\in N$ if and only if $u\in V$ and
\begin{equation*}
\int_{\partial\Omega_{1}}t_1(u^1)\cdot v^1\,{\rm d}\mathcal{H}^{d-1}+\int_{\partial\Omega_{2}}t_2(u^2)\cdot v^2\,{\rm d}\mathcal{H}^{d-1}=0
\end{equation*}
holds true for every $v\in V$. The last condition is equivalent to
\begin{equation}\label{co1}
\int_{S_1}t_1(u^1)\cdot v^1\,{\rm d}\mathcal{H}^{d-1}=0,
\end{equation}
\begin{equation}\label{co2}
\int_{S_2}t_2(u^2)\cdot v^2\,{\rm d}\mathcal{H}^{d-1}=0,	
\end{equation}
and
\begin{equation}\label{co3}
\int_\Gamma t_1(u^1)\cdot v^1+t_2(u^2)\cdot v^2\,{\rm d}\mathcal{H}^{d-1}=0.
\end{equation}
Since $T(v^1)=T(v^2)$, i.e. $v^1=v^2$ on $\Gamma$, \eqref{co1}, \eqref{co2} and \eqref{co3} are equivalent to
\begin{equation*}
\begin{cases}
u^1=u^2&\hspace{.1in}\mbox{on}\hspace{.1in}\Gamma\\
t_1(u^1)+t_2(u^2)=0&\hspace{.1in}\mbox{on}\hspace{.1in}\Gamma
\end{cases}
\end{equation*}
and
\begin{equation*}
\begin{cases}
t_1(u^1)=0&\hspace{.1in}\mbox{on}\hspace{.1in}S_1\\
t_2(u^2)=0&\hspace{.1in}\mbox{on}\hspace{.1in}S_2.
\end{cases}
\end{equation*}
This analysis completes the proof. The next section focuses on some examples of fixed transmission problems for which a regularity theory is available.

\subsection{Some developments in regularity theory}\label{subsec_tar}

Important advances in the regularity theory of fixed transmission problems are in \cite{Li-Vogelius2000}. In that paper, the authors derive global $W^{1,\infty}$-estimates and $C^{1,\alpha}$-regularity for solutions to divergence form elliptic equations with piecewise H\"older continuous coefficients. We proceed by introducing some notation and making the assumptions more precise. 

Let $\Omega\subset\mathbb{R}^d$ be a bounded domain with $C^{1,\alpha}$-regular boundary, $\alpha\in(0,1)$, and let $\Omega_m$, $1\le m\le L$, be a finite number of disjoint subdomains of $\Omega$, each with $C^{1,\alpha}$ boundary. Furthermore, suppose that 
\[
    \overline{\Omega}=\bigcup_{m=1}^L\overline{\Omega_m}.
\]
Given $\overline{x}\in\overline{\Omega}$, one supposes that there exists $r>0$ and an appropriate rotation of our fixed coordinate system, such that the set
\[
   \left (\bigcup_{m=1}^L\partial\Omega_m\right)\cap B_r(\overline{x})
\]
consists of the graphs of a finite number of $C^{1,\alpha}$-regular functions. Denote with $\ell(\overline{x},r)$ the number of these functions, and let $K(\overline{x},r)$ denote the maximum of their $C^{1,\alpha}$ norms. The quantity
\begin{equation*}
\mathcal{K}:=\sup_{\overline{x}\in\Omega}\inf_{r>0}\left(K(\overline{x},r)+\ell(\overline{x},r)+\frac{1}{r}\right)
\end{equation*}
is referred to as the $C^{1,\alpha}$ modulus of the total boundary set $(\cup_{m=1}^L\partial\Omega_m)$.

For $\mu\in(0,1)$, let $A^{(m)}\in C^\mu\left(\overline{\Omega}_m\right)$  be a symmetric, positive definite matrix-valued function, and define
\begin{equation}\label{assum1}
A(x):=A^{(m)}(x)\hspace{.1in}\mbox{in}\hspace{.1in}x\in\Omega_m,\,1\le m\le L.
\end{equation}
Denote with $0<\overline{\lambda}<\overline{\Lambda}<\infty$ the ellipticity constants associated with $A$. Similarly, let $g^{(m)}\in C^{0,\mu}\left(\overline{\Omega};\mathbb{R}^d\right)$, and define
\begin{equation}\label{assum2}
g(x):=g^{(m)}(x)\hspace{.1in}\mbox{in}\hspace{.1in}x\in\Omega_m,\,1\le m\le L.
\end{equation}
Finally, suppose
\begin{equation}\label{assum3}
h\in L^\infty(\Omega)
\end{equation}
and
\begin{equation}\label{assum4}
\varphi\in C^{1,\mu}(\partial\Omega).
\end{equation}
The first of the results in \cite{Li-Vogelius2000} concerns $C^{1,\alpha'}$ interior estimates. In brief, this result asserts that the restriction of the solution $u$ to each subdomain $\Omega_m$ can be extended to $\Omega_\varepsilon:=\{x\in\Omega:\,{\rm dist}(x,\partial\Omega)>\varepsilon\}$ as a $C^{1,\alpha'}$ function, with a norm that is independent of the distances between the subdomain interfaces.

\begin{Theorem}\label{interior}
Suppose \eqref{assum1}-\eqref{assum3} hold. Suppose $\alpha'$ satisfies $0<\alpha'\le\mu$ and $\alpha'<\alpha/((\alpha+1)d)$, and let $\varepsilon>0$. Then there exists $C>0$, depending only on $\Omega,d,\alpha,\alpha',\varepsilon,\overline{\lambda},\overline{\Lambda},\left\|A^{(m)}\right\|_{C^{0,\alpha'}(\overline{\Omega})}$, and $\mathcal{K}$, such that, if $u\in W^{1,2}(\Omega)$ is a solution to
\begin{equation*}
\partial_i(A_{ij}\partial_ju)=h+\partial_ig_i\hspace{.1in}\mbox{in}\hspace{.1in}x\in\Omega,
\end{equation*}
then
\begin{align*}
\max_{1\le m\le L}\|u\|_{C^{1,\alpha'}(\overline{\Omega}_m\cap\Omega_\varepsilon)}\le&C\Big(\|u\|_{L^\infty(\Omega)}+\|h\|_{L^\infty(\Omega)}\\
&+\max_{1\le m\le L}\left\|g^{(m)}\right\|_{C^{0,\alpha'}(\overline{\Omega})}\Big).
\end{align*}
\end{Theorem}
The variant of Theorem \ref{interior} covering boundary estimates is the following.
\begin{Theorem}
Suppose \eqref{assum1}-\eqref{assum4} hold. Suppose $\alpha'$ satisfies $0<\alpha'\le\mu$ and $\alpha'<\alpha/((\alpha+1)d)$, and let $\varepsilon,r>0$. Then there exists $C>0$, depending only on $\Omega,d,\alpha,\alpha',\varepsilon,\overline{\lambda},\overline{\Lambda},\left\|A^{(m)}\right\|_{C^{0,\alpha'}(\overline{\Omega})}$ and $\mathcal{K}$, such that, if for some $\overline{x}\in\partial{\Omega}$, $u\in W^{1,2}(\Omega\cap B_{2r}(\overline{x}))$ is a solution to
\begin{equation*}
\begin{cases}
\partial_i(A_{ij}\partial_ju)=h+\partial_ig_i&\hspace{.1in}\mbox{in}\hspace{.1in}x\in\Omega\cap B_{2r}(\overline{x})\\
u=\varphi&\hspace{.1in}\mbox{on}\hspace{.1in}x\in\partial\Omega\cap B_{2r}(\overline{x}),
\end{cases}
\end{equation*}
then
\begin{align*}
\max_{1\le m\le L}\|u\|_{C^{1,\alpha'}(\overline{\Omega}_m\cap B_r(\overline{x}))}&\le C\Big(\|u\|_{L^\infty(\Omega\cap B_r(\overline{x}))}+\|\varphi\|_{C^{1,\alpha'}(\partial\Omega\cap B_{2r}(\overline{x}))}\\
&\quad+\|h\|_{L^\infty(\Omega\cap B_r(\overline{x}))}\\
&\quad+\max_{1\le m\le L}\left\|g^{(m)}\right\|_{C^{0,\alpha'}(\overline{\Omega}\cap B_r(\overline{x}))}\Big).
\end{align*}
\end{Theorem}
Combining the above interior and boundary estimates with the maximum principle (see \cite{Gilbarg-Trudinger}) we arrive at the following $C^{1,\alpha'}$ global estimate.

\begin{Corollary}
Suppose \eqref{assum1}-\eqref{assum4} hold. Suppose $\alpha'$ satisfies $0<\alpha'\le\mu$ and $\alpha'<\alpha/((\alpha+1)d)$. Then there exists $C>0$, depending only on $\Omega,d,\alpha,\alpha',\lambda,\overline{\Lambda},\left\|A^{(m)}\right\|_{C^{0,\alpha'}(\overline{\Omega})}$, and $\mathcal{K}$, such that, if $u\in W^{1,2}(\Omega)$ is a solution to
\begin{equation*}
\begin{cases}
\partial_i(A_{ij}\partial_ju)=h+\partial_ig_i&\hspace{.1in}\mbox{in}\hspace{.1in}x\in\Omega\\
u=\varphi&\hspace{.1in}\mbox{on}\hspace{.1in}x\in\partial\Omega,
\end{cases}
\end{equation*}
then
\begin{align*}
\max_{1\le m\le L}\|u\|_{C^{1,\alpha'}(\overline{\Omega}_m)}&\le C\Big(\|\varphi\|_{C^{1,\alpha'}(\partial\Omega}+\|h\|_{L^\infty(\Omega)}\\
&\quad+\max_{1\le m\le L}\left\|g^{(m)}\right\|_{C^{0,\alpha'}(\overline{\Omega})}\Big).
\end{align*}
\end{Corollary}

We now describe the methods of the proof in \cite{Li-Vogelius2000}. To keep it brief and as clear as possible, we restrict our attention to Theorem \ref{interior} in the case $h=0$ and $g=0$. Let $\mathcal{A}\left(\overline{\lambda},\overline{\Lambda}\right)$ denote the set of measurable, symmetric, positive definite matrix-functions $A(x)$ satisfying
\begin{equation*}
\overline{\lambda}Id\le A(x)\le\overline{\Lambda}Id.
\end{equation*}
Define a scaling invariant subclass of $\mathcal{A}\left(\overline{\lambda},\overline{\Lambda}\right)$, denoted by $\overline{\mathcal{A}}\left(\overline{\lambda},\overline{\Lambda}\right)$, as follows. Fix $\ell\in\mathbb{N}$ and denote with $\{L_1,\dots,L_\ell\}$ a collection of $\ell$ parallel hyperplanes in $\mathbb{R}^d$, dividing $\mathbb{R}^d$ in $\ell+1$ regions. Denote such regions with $R_1,\dots,R_{\ell+1}$. Let $\overline{A}^{(1)},\dots,\overline{A}^{(\ell+1)}$ be any $\ell+1$ symmetric, positive definite constant matrices in $\mathcal{A}\left(\overline{\lambda},\overline{\Lambda}\right)$ and define
\begin{equation*}
\overline{A}_{ij}(x):=\overline{A}_{ij}^m\hspace{.1in}\mbox{on}\hspace{.1in}x\in R_m,\,1\le m\le l+1.
\end{equation*}
The subclass $\overline{\mathcal{A}}\left(\overline{\lambda},\overline{\Lambda}\right)$ consists of all such matrix functions $\overline{A}$.

The classical Schauder estimates, Cordes-Nirenberg estimates, and $W^{2,p}$ estimates can be viewed as perturbation theories from the corresponding estimates for solutions to the Laplace equation. Here, the approach is analogous, as the problem of interest is regarded  as a perturbation of
\begin{equation}\label{help}
\partial_i\left(\overline{A}_{ij}(x)\partial_jv\right)=0,
\end{equation}
with $\overline{A}\in\overline{\mathcal{A}}\left(\overline{\lambda},\overline{\Lambda}\right)$.

To establish Theorem \ref{interior} in the case $h=0$ and $g=0$, one first studies elliptic regularity estimates for solutions to \eqref{help}. Although the hyperplanes in the definition of $\overline{\mathcal{A}}\left(\overline{\lambda},\overline{\Lambda}\right)$ are allowed to get arbitrarily close to each other, estimates that are uniform in $\overline{A}\in\overline{\mathcal{A}}\left(\overline{\lambda},\overline{\Lambda}\right)$ are available. Hence, the Caccioppoli inequality and the interior De Giorgi-Nash estimates yield bounds for all derivatives of the solution $v$ in each region $\overline{R}_m$. More precisely, with $\Omega=(-1,1)^d$, one concludes that, for any positive integer $k$, any $\varepsilon>0$, any $\overline{A}\in\overline{\mathcal{A}}\left(\overline{\lambda},\overline{\Lambda}\right)$ and any solution $v$ to \eqref{help}, it holds
\begin{equation}\label{starting}
\max_{1\le m\le l+1}\|v\|_{C^k(\overline{R}_m\cap(1-\varepsilon)\Omega)}\le C\|v\|_{L^\infty(\Omega)}.
\end{equation}
Starting from \eqref{starting}, perturbation methods allow us to show that, for any $q>d$ and $\overline{\alpha}\in(0,1)$, there exists $\varepsilon_0>0$, depending only on $d,q,\overline{\alpha},\overline{\lambda}$ and $\overline{\Lambda}$, such that if $A\in\mathcal{A}\left(\overline{\lambda},\overline{\Lambda}\right)$ and $\overline{A}\in\overline{\mathcal{A}}\left(\overline{\lambda},\overline{\Lambda}\right)$ satisfy
\begin{equation}\label{closeness}
\sup_{r\in(0,1)}r^{-\overline{\alpha}}\left(\frac{1}{|r\Omega|}\int_{r\Omega}\left|A(x)-\overline{A}(x)\right|^q\,{\rm d}x\right)^{\frac{1}{q}}\le\varepsilon_0
\end{equation}
and $u\in W^{1,2}(\Omega)$ is a solution to
\begin{equation*}
\partial_i(A_{ij}(x)\partial_{j})u=0\hspace{.1in}\mbox{in}\hspace{.1in}\Omega,
\end{equation*}
with
\begin{equation*}
\|u\|_{L^\infty(\Omega)}\le1,
\end{equation*}
then there exists a continuous piecewise linear solution $p$ to
\begin{equation*}
\partial_i\left(\overline{A}_{ij}(x)\partial_jp\right)=0\hspace{.1in}\mbox{in}\hspace{.1in}\frac{1}{4}\Omega,
\end{equation*}
whose coefficients are bounded in absolute value by a constant $C$, depending only on $d,q,\overline{\alpha},\overline{\lambda}$, and $\overline{\Lambda}$, such that
\begin{equation*}
|u(x)-p(x)|\le C|x|^{1+\overline{\alpha}}\hspace{.1in}\mbox{in}\hspace{.1in}\frac{1}{4}\Omega.
\end{equation*}

Under the hypotheses of Theorem \ref{interior}, the condition \eqref{closeness} can be verified at every point $\overline{x}\in\Omega$ by translation and dilation, and by appropriate selection of $\overline{A}_{\overline{x}}\in\overline{\mathcal{A}}\left(\overline{\lambda},\overline{\Lambda}\right)$. The $L^\infty$-interior estimates for the gradient of solutions to the equation $\partial_i(A_{ij}\partial_ju)=0$ follow immediately. The H\"older interior estimates for the gradient require some further work, since at different points $\overline{x}\in\Omega$, the orientation of the hyperplanes associated with the matrices $\overline{A}_{\overline{x}}(x)$ differ by a rotation, determined by the geometry of $\Omega_m$, and since $p_{\overline{x}}$ is only piecewise linear given a fixed set of planes.

\subsection{Geometry of the interface and optimal regularity}\label{subsec_tango}

So far we have dealt with transmission problems with smooth fixed interfaces. One of the main novelties presented in \cite{CSCS2021} is that the interface has only $C^{1,\alpha}$-regularity. In what follows, we detail the setting of that paper and describe its main contributions. 

Let $\Omega\subseteq\mathbb{R}^d$ be a bounded domain. Fix $\Omega_1\Subset\Omega$, and suppose that $\Gamma:=\partial\Omega_1$ is a $C^{1,\alpha}$ manifold for some $\alpha\in(0,1)$. Set $\Omega_2:=\Omega\setminus\Omega_1$. The authors in that paper consider the problem of finding a continuous function $u:\overline{\Omega}\to\mathbb{R}$ such that
\begin{equation}\label{1}
\begin{cases}
\Delta u_1=0&\hspace{.1in}\mbox{in}\hspace{.1in}\Omega_1\\
\Delta u_2=0&\hspace{.1in}\mbox{in}\hspace{.1in}\Omega_2\\
u_2=0&\hspace{.1in}\mbox{in}\hspace{.1in}\partial\Omega\\
u_1=u_2&\hspace{.1in}\mbox{on}\hspace{.1in}\Gamma\\
(u_1)_\nu-(u_2)_\nu=g&\hspace{.1in}\mbox{on}\hspace{.1in}\Gamma
\end{cases}
\end{equation}
where $u_1:=u\big|_{\overline{\Omega_1}}$, $u_2:=u\big|_{\overline{\Omega_2}}$, $g\in C^{0,\alpha}(\Gamma)$, and $\nu$ is the unit normal vector on $\Gamma$ that is interior to $\Omega_1$. The last two equations in \eqref{1} are called \textit{transmission conditions}. Note that if $g\equiv0$ then $u$ is a harmonic function in $\Omega$. Hence, in order to have a meaningful transmission condition, we suppose that
\begin{equation*}
g(x)\ge0\hspace{.1in}\mbox{for}\hspace{.1in}x\in\Gamma.
\end{equation*}
The main result in \cite{CSCS2021} is the following theorem.
\begin{Theorem}\label{21}
There exists a unique solution $u$ to \eqref{1}. Moreover, $u_1\in C^{1,\alpha}\left(\overline{\Omega_1}\right)$, $u_2\in C^{1,\alpha}\left(\overline{\Omega_2}\right)$, and there exists $C=C(d,\alpha,\Gamma)>0$ such that
\begin{equation*}
\|u_1\|_{C^{1,\alpha}(\overline{\Omega_1})}+\|u_2\|_{C^{1,\alpha}(\overline{\Omega_2})}\le C\|g\|_{C^{0,\alpha}(\Gamma)}.
\end{equation*}
In addition, $u\in C^{0,{\rm Log-Lip}}\left(\overline{\Omega}\right)$.
\end{Theorem}
The existence and the uniqueness easily follow from a representation formula of the solution via the Green function $G$ associated with the Laplace operator in $\Omega$
\begin{equation}\label{3}
u(x)=\int_{\Gamma}G(x,y)g(y)\,{\rm d}\mathcal{H}^{d-1}\hspace{.1in}\mbox{for}\hspace{.1in}x\in\Omega.
\end{equation}
Thanks to the well-known properties of the Green function we achieve, using the above formula, that $u\in C^{0,{\rm Log-Lip}}\left(\overline{\Omega}\right)$. The main issue is the regularity of $u_i$ up to the interface for $i=1,2$. To establish this fact one resorts to a number of building blocks. We continue by describing them and start with a discussion about the \emph{flat problem}.

For $a\in\mathbb{R}$, we denote
\begin{align*}
x=&(x',x_d)\in\mathbb{R}^{d-1}\times\mathbb{R}\\
B_{r,a}=&B_r(0',a)\\
B_{r,a}^+=&B_r(0',a)\cap\{x_d>a\}\\
B_{r,a}^-=&B_r(0',a)\cap\{x_d<a\}\\
T_{r,a}=&\{x\in B_r(0',a):\,x_d=a\}\\
T_a=&B_1\cap\{x_d=a\}\\
T_a^+=&B_1\cap\{x_d\ge a\}\\
T_a^-=&B_1\cap\{x_d\le a\}.
\end{align*}
When $a=0$, we write $T=T_0$ and $B_r^{\pm}=B_{r,0}^{\pm}$. In the context of flat interfaces, the analogous of Theorem \ref{21} reads as follows.
\begin{Theorem}\label{2}
Let $r>0$ and $a\in\mathbb{R}$. Given $\alpha,\gamma\in(0,1)$, let $g\in C^{0,\alpha}(T_{r,a})$ and $f\in C^{0,\gamma}\left(\overline{B_{r,a}}\right)$. Then there exists a unique solution $v\in C^\infty(B_{r,a}\setminus T_{r,a})\cap C^{0,\gamma}\left(\overline{B_{r,a}}\right)$ to the flat transmission problem
\[
\begin{cases}
\Delta v=g\,{\rm d}\mathcal{H}^{d-1}\big|_{T_{r,a}}&\hspace{.1in}\mbox{in}\hspace{.1in}B_{r,a}\\
v=f&\hspace{.1in}\mbox{in}\hspace{.1in}\partial B_{r,a}.
\end{cases}
\]
Moreover, if $v^\pm=v\chi_{\overline{B_{r,a}^\pm}}$, then $v\in C^{1,\alpha}\left(\overline{B_{r/2,a}^\pm}\right)$ and
\begin{equation*}
\left\|v^\pm\right\|_{C^{1,\alpha}(\overline{B_{r/2,a}^\pm})}\le C\left(\|g\|_{C^{0,\alpha}(T_{r,a})}+\|f\|_{L^\infty(\partial B_{r,a})}\right),
\end{equation*}
where $C$ depends only on $d$, $\alpha$, and $r$. If $g\in C^{k-1,\alpha}(T_{r,a})$, $k\ge1$, then $v\in C^{k,\alpha}\left(\overline{B_{r/2,a}^\pm}\right)$ and
\begin{equation*}
\left\|v^\pm\right\|_{C^{k,\alpha}(\overline{B_{r/2,a}^\pm})}\le C\left(\|g\|_{C^{k-1,\alpha}(T_{r,a})}+\|f\|_{L^\infty(\partial B_{r,a})}\right).
\end{equation*}
\end{Theorem}
Briefly, the proof of the above theorem relies on the construction of two Dirichlet-Neumann problems, one for each hemisphere. Then, the solution to the original problem is given by the sum of the solutions of the aforementioned two problems.

\begin{Corollary}\label{6}
Given $|a|<1/4$, $c_0>0$ and $f\in C^{0,\gamma}\left(\overline{B_{1}}\right)$, with $\gamma\in(0,1)$, there exists a unique solution $v\in C^\infty(B_1\setminus T_a)\cap C^{0,\gamma}\left(\overline{B_{1}}\right)$ to
\[
\begin{cases}
\Delta v=c_0\,{\rm d}\mathcal{H}^{d-1}&\hspace{.1in}\mbox{in}\hspace{.1in}B_1\\
v=f&\hspace{.1in}\mbox{on}\hspace{.1in}\partial B_1
\end{cases}
\]
such that, for any $k\ge1$,
\begin{equation*}
\left\|v^\pm\right\|_{C^{k,\alpha}(\overline{B_{1/2}}\cap T_a^{\pm})}\le C\left(c_0+\|f\|_{L^\infty(\partial B_1)}\right),
\end{equation*}
where $C=C(d,\alpha,k)>0$.
\end{Corollary}

Before continuing, further notation is required. Fix $\varepsilon>0$, and let $\Omega_\varepsilon=\{x\in\Omega:\,\textnormal{dist}(x,\partial\Omega)<\varepsilon\}$ and $\Gamma_\varepsilon=\{x\in\Omega:\,\textnormal{dist}(x,\Gamma)<\varepsilon\}$. Consider the average
\begin{equation*}
u_\varepsilon(x)=\frac{1}{|B_\varepsilon(x)|}\int_{B_\varepsilon(x)}u\,{\rm d}y\hspace{.1in}\mbox{for}\hspace{.1in}x\in\Omega_\varepsilon.
\end{equation*}

\begin{Proposition}[Properties of averages]\label{averages}
Let $u$ be the distributional solution given by \eqref{3}. The following properties hold
\begin{enumerate}
\item[$(i)$] if $B_\varepsilon(x)\cap\Gamma=\emptyset$ then $u_\varepsilon(x)=u(x)$;
\item[$(ii)$] $u_\varepsilon\to u$ locally uniformly in $\Omega$ as $\varepsilon\to0^+$;
\item[$(iii)$] if $g\in L^\infty(\Omega)$ then $g_\varepsilon\in C_c(\Gamma_\varepsilon)$, where
\begin{equation*}
g_\varepsilon(x)=\frac{1}{|B_\varepsilon|}\int_{\Gamma\cap B_\varepsilon(x)}g\,{\rm d}\mathcal{H}^{d-1}\hspace{.1in}\mbox{for}\hspace{.1in}x\in\Gamma_\varepsilon.
\end{equation*}
Moreover, $\Delta u_\varepsilon(x)=g_\varepsilon(x)$, for every $x\in\Omega_\varepsilon$.
\end{enumerate}
\end{Proposition}
The proof of the above proposition is quite simple, and it is essentially based on well-known properties of harmonic functions, the Lebesgue Dominated Convergence Theorem, and a change of variables. The following result is instrumental in the analysis.
\begin{Lemma}
Let $\Gamma$ be as in Theorem \ref{stability}. Define $M:=1+2\theta$ and let $x\in B_{1-M\varepsilon}$ be such that $\textnormal{dist}(x,\Gamma)<\varepsilon$. Then
\begin{align}\label{4}
\left\{y':\,\big(y',\psi(y')\big)\in B_\varepsilon(x)\right\}\subseteq&B_{((M\varepsilon)^2-(x_d+\theta\varepsilon)^2)^{1/2}}'(x') \notag\\
=&\{y':\,(y',-\theta\varepsilon)\in B_{M\varepsilon}(x)\}
\end{align}
and
\begin{align}\label{5}
\left\{y':\,\big(y',\psi(y')\big)\in B_{M\varepsilon}(x)\right\}\supseteq&B_{(\varepsilon^2-(x_d+\theta\varepsilon)^2)^{1/2}}'(x') \notag\\
=&\{y':\,(y',-\theta\varepsilon)\in B_{\varepsilon}(x)\}.
\end{align}
\end{Lemma}

\begin{Theorem}[Stability]\label{stability}
Let $\varepsilon >0,\theta<1/2$ and $0<\delta,\gamma<1$ be given. Let $\Gamma=\left\{\big(y',\psi(y')\big):\,y'\in B_1'\right\}$, where $\psi$ is a Lipschitz function. Suppose $\Gamma$ is $\theta\varepsilon$-flat in $B_1$, in the sense that
\begin{equation*}
\Gamma\subseteq\{x\in B_1:\,|x_d|<\theta\varepsilon\},
\end{equation*}
and that $\Gamma$ is also $\varepsilon$-horizontal in $B_1$. That is
\begin{equation*}
1-\varepsilon\le\nu(x)\cdot(0',1)=\left|1+|D'\psi(x')|^2\right|^{-\frac{1}{2}}\le1,
\end{equation*}
for every $x\in\Gamma$, where $\nu(x)$ denotes the upward pointing normal on $\Gamma$. Then there exists $C=C(d)>0$ such that, for any $u\in C^{0,\gamma}\left(\overline{B_{1}}\right)$ and $g\in L^\infty(\Gamma)$ satisfying
\[
\begin{cases}
\Delta u=g\,{\rm d}\mathcal{H}^{d-1}\big|_{\Gamma}&\hspace{.1in}\mbox{in}\hspace{.1in}B_1\\
|g-1|\le\delta&\hspace{.1in}\mbox{on}\hspace{.1in}\Gamma,
\end{cases}
\]
the classical solution $v\in C^\infty(B_1\setminus T_{-\theta\varepsilon})\cap C^{0,\gamma}\left(\overline{B_{1}}\right)$ to the flat problem
\[
\begin{cases}
\Delta v={\rm d}\mathcal{H}^{d-1}\big|_{T_{-\theta\varepsilon}}&\hspace{.1in}\mbox{in}\hspace{.1in}B_1\\
v=u&\hspace{.1in}\mbox{on}\hspace{.1in}\partial B_1
\end{cases}
\]
satisfies
\begin{equation*}
|u-v|\le C\left(\theta+\delta+\varepsilon^\gamma\right)\quad\textnormal{in }B_{1/2}.
\end{equation*}
\end{Theorem}

\begin{Remark}\label{13}
The interface for the flat problem in Theorem \ref{stability} is $T_{-\theta\varepsilon}=B_1\cap\{x_d=-\theta\varepsilon\}$, which lies below $\Gamma$ in the $x_d$-direction. To approximate $u$ with the solution to a flat problem where the interface lies above $\Gamma$ in the $x_d$-direction, it is enough to consider the classical solution $v$ to
\[
\begin{cases}
\Delta v={\rm d}\mathcal{H}^{d-1}\big|_{T_{\theta\varepsilon}}&\hspace{.1in}\mbox{in}\hspace{.1in}B_1\\
v=u&\hspace{.1in}\mbox{on}\hspace{.1in}\partial B_1.
\end{cases}
\]
In this case, the same conclusion as in Theorem \ref{stability} holds.
\end{Remark}
The proof of the stability result is quite technical, and it is obtained via a novel geometric approach which is based on the mean value property and the maximum principle. What the stability tells us is that if the flatness and oscillation of the interface $\Gamma$ are controlled, then one can construct a solution for a flat interface problem, where the flat interface does not intersect $\Gamma$. Also, it is possible to quantify how close solutions are, depending only on the geometric properties of $\Gamma$ and the basic regularity of $u$.

Now, suppose that $\Gamma$ is an interface in $B_1$ given by the graph of a function $x_d=\psi(x'):T\to\mathbb{R}$. Thus, we can write $B_1=\Omega_1\cup\Gamma\cup\Omega_2$, where $\Omega_1=\{x=(x',x_d)\in B_1:\,x_d>\psi(x')\}$. It is also natural to suppose that $0\in\Gamma$.

\begin{Lemma}\label{15}
Given $0<\alpha,\gamma<1$, there exist $C_0>0$, $\lambda\in(0,1/2)$, $0<\theta,\delta,\varepsilon<\lambda$, depending only on $d,\alpha$ and $\gamma$, such that, for every $u\in C^{0,\gamma}\left(\overline{B_{1}}\right)$ satisfying
\[
\begin{cases}
\Delta u=g\,{\rm d}\mathcal{H}^{d-1}\big|_{\Gamma}&\hspace{.1in}\mbox{in}\hspace{.1in}B_1\\
|u|\le1&\hspace{.1in}\mbox{in}\hspace{.1in}B_1\\
|g-1|\le\delta&\hspace{.1in}\mbox{on}\hspace{.1in}\Gamma,
\end{cases}
\]
if $\Gamma$ is $\theta\varepsilon$-flat and $\varepsilon$-horizontal in $B_1$, then there are linear polynomials $P_1(x)=A\cdot x+B$ and $Q_1(x)=C\cdot x+B$, with $A,C\in\mathbb{R}^d$, $B\in\mathbb{R}$ and $|A|+|B|+|C|\le C_0$, for which
\[
|u_1(x)-P_1(x)|\le\lambda^{1+\alpha}\hspace{.1in}\mbox{for}\hspace{.1in}x\in\Omega_1\cap B_\lambda
\]
and
\[
|u_2(x)-Q_1(x)|\le\lambda^{1+\alpha}\hspace{.1in}\mbox{for}\hspace{.1in}x\in\Omega_2\cap B_\lambda.
\]
Moreover, $D'P_1=D'Q_1'$ and $(P_1)_{x_d}-(Q_1)_{x_d}=1$.
\end{Lemma}

\begin{Lemma}\label{16}
Given $\alpha\in(0,1)$, there exist $C_0>0$, $\lambda\in(0,1/2)$, $\delta\in(0,1)$ depending only on $d$ and $\alpha$, such that, for any distributional solution $u\in C\left(\overline{B_{1}}\right)$ to
\[
\begin{cases}
\Delta u=g\,{\rm d}\mathcal{H}^{d-1}\big|_{\Gamma}&\hspace{.1in}\mbox{in}\hspace{.1in}B_1\\
|u|\le1&\hspace{.1in}\mbox{in}\hspace{.1in}B_1\\
|g|\le\delta&\hspace{.1in}\mbox{on}\hspace{.1in}\Gamma,
\end{cases}
\]
there is a linear polynomial $P(x)=A\cdot x+B$, with $A\in\mathbb{R}^d$ and $B\in\mathbb{R}$ and $|A|+|B|\le C_0$, satisfying
\begin{equation*}
|u(x)-P(x)|\le\lambda^{1+\alpha}\hspace{.1in}\mbox{for}\hspace{.1in}x\in B_\lambda.
\end{equation*}
\end{Lemma}
The previous two lemmas are key tools in the following theorem. In fact, the proof of the next result is based on a clever induction argument that involves a sequence of scaled fixed transmission problems. The very general idea is that flat solutions are asymptotically close to non-flat solutions. The theorem reads as follows.

\begin{Theorem}[Pointwise $C^{1,\alpha}$ boundary regularity]\label{22}
Let 
\[
    \Gamma=\left\{\big(y',\psi(y')\big):\,y'\in B_1'\right\},
\]
where $\psi$ is a $C^{1,\alpha}$-regular function, for some $\alpha\in(0,1)$. Suppose that $0\in\Gamma$. Let $u$ be a distributional solution to the transmission problem
\begin{equation*}
\Delta u=g\,{\rm d}\mathcal{H}^{d-1}\big|_{\Gamma},
\end{equation*}
where $g\in L^\infty(\Gamma)\cap C^{0,\alpha}(0)$ is nonnegative. Then there are linear polynomials $P(x)=A\cdot x+B$ and $Q(x)=C\cdot x+B$ such that
\begin{align*}
|u_1(x)-P(x)|\le&D|x|^{1+\alpha}\hspace{.1in}\mbox{for}\hspace{.1in}x\in\Omega_1\cap B_{1/2},\\
|u_2(x)-Q(x)|\le&D|x|^{1+\alpha}\hspace{.1in}\mbox{for}\hspace{.1in}x\in\Omega_2\cap B_{1/2},
\end{align*}
with
\begin{equation*}
|A|+|B|+|C|+D\le C_0\|\psi\|_{C^{1,\alpha}(B_1')}\left([g]_{C^{0,\alpha}(0)}+\|g\|_{L^\infty}(\Gamma)\right)
\end{equation*}
and $C_0=C_0(d,\alpha)>0$.
\end{Theorem}

To prove Theorem \ref{21} one resorts to Campanato’s characterization of $C^{1,\alpha}$-spaces and a technical result that patches the interior and boundary estimates together.

\begin{Theorem}[Campanato's characterization of $C^{1,\alpha}$-spaces]\label{23}
Let $u$ be a measurable function defined on a bounded $C^{1,\alpha}$ domain $\Omega$. Then $u\in C^{1,\alpha}\left(\overline{\Omega}\right)$ if and only if there exists $C_0>0$ such that, for every $x\in\overline{\Omega}$, one finds a linear polynomial $Q_x(z)$ satisfying
\begin{equation*}
|u(z)-Q_x(z)|\le C_0|z-x|^{1+\alpha},
\end{equation*}
for every $z\in B_1(x)\cap\Omega$. In this case, if $C_*$ denotes the least constant $C_0$ for which the property above holds, we have
\begin{equation*}
\|u\|_{C^{1,\alpha}(\overline{\Omega})}\sim C_*+\sup_{x\in\overline{\Omega}}|Q_x|,
\end{equation*}
where $|Q_x|$ denotes the sum of the coefficients of the polynomial $Q_x(z)$.
\end{Theorem}

At this point, the proof of Theorem \ref{21} proceeds as follows. We first introduce three conditions and relate them with the conclusion of Theorem \ref{21} through a proposition. In fact, the proposition ensures that if $u_1$ and $u_2$ satisfy the aforementioned conditions then the regularity up to the fixed interface is available. We continue by introducing the conditions of interest.

From now on, let $S$ be a collection of measurable functions defined on a bounded $C^{1,\alpha}$-regular domain $U\subset\mathbb{R}^d$.

\begin{Condition}[Interior estimates]\label{cond_1}
Let $u\in S$ and $d_x:=\textnormal{dist}(x,\partial\Omega)$. There exist $A,C,D>0$ such that, for every $x\in U$, one can find a linear polynomial $P_x(z)$ satisfying
\begin{equation*}
\|P_x\|_{L^\infty(B)}+d_x\|DP_x\|_{L^\infty(B)}\le C\|u\|_{L^\infty(B)},
\end{equation*}
and
\begin{equation*}
|u(z)-P_x(z)|\le\left(A\frac{\|u\|_{L^\infty(B)}}{d_x^{1+\alpha}}+D\right)|z-x|^{1+\alpha},
\end{equation*}
for every $z\in B:=B_{d_x/2}(x)\subset U$.
\end{Condition}

\begin{Condition}[Boundary estimates]\label{cond_2}Let $u\in S$. There exists $E>0$ such that, for every $y\in\partial U$, there is a linear polynomial $P_y(z)$ satisfying
\begin{equation*}
\|P_y\|_{L^\infty(U)}+\|DP_y\|_{L^\infty(U)}\le E
\end{equation*}
and
\begin{equation*}
|u(z)-P_y(z)|\le E|z-y|^{1+\alpha},
\end{equation*}
for every $z\in\overline{U}$.
\end{Condition}

\begin{Condition}[Invariance property]\label{cond_3} For every $u\in S$, and for every $y\in\partial U$ with a corresponding linear polynomial $P_y$ as in Condition \ref{cond_2}, the function $v=u-P_y$ also satisfies the estimates in Condition \ref{cond_1}.
\end{Condition}

\begin{Proposition}\label{24}
Let $S$ be a collection of measurable functions defined on a bounded $C^{1,\alpha}$-regular domain $U$. For $x\in U$, we let $d_x=\textnormal{dist}(x,\partial U)$. Fix $u\in S$, and suppose Conditions \ref{cond_1}-\ref{cond_3} are in force. Then $S\subseteq C^{1,\alpha}\left(\overline{U}\right)$, and there exists $C>0$, depending only on $A,C,D,E$, such that
\begin{equation*}
\|u\|_{C^{1,\alpha}(\overline{U})}\le C\|u\|_{L^\infty(U)}.
\end{equation*}
\end{Proposition}

Now, we show that $u_1$ and $u_2$ verify Conditions \ref{cond_1}-\ref{cond_3}. This fact builds upon Proposition \ref{24} to produce the conclusion of Theorem \ref{21}, completing its proof.

Let  $u\in C^{0,{\rm Log-Lip}}\left(\overline{\Omega}\right)$ be given by \eqref{3}. We consider only $u_2:\overline{\Omega_2}\to\mathbb{R}$, as a similar argument yields the result for $u_1:\overline{\Omega_1}\to\mathbb{R}$ and start with Condition \ref{cond_1}.

Fix $x\in\Omega_{2}$. Since $u_2$ is harmonic, it is smooth in $\Omega_{2}$, so we can define
\begin{equation*}
P_x(z):=u_2(x)+Du_2(x)\cdot(z-x).
\end{equation*}
Then, by classical interior estimates for harmonic functions,
\begin{align*}
\|P_x\|_{L^\infty(B)}+d_x\|DP_x\|_{L^\infty(B)}\le&\|u_2\|_{L^\infty(B)}+d_x\|Du_2\|_{L^\infty(B)}\\
&+d_x\|Du_2\|_{L^\infty(B)}\\
\le&\|u_2\|_{L^\infty(B)}+2d\|u_2\|_{L^\infty(B)}\\
\le&(1+2d)\|Du_2\|_{L^\infty(B)}.
\end{align*}
Moreover,
\begin{align*}
|u_2(z)-P_x(z)|\le&\left\|D^2u_2\right\|_{L^\infty(B)}|z-x|^2\\
\le&d\frac{\|u_2\|_{L^\infty(B)}}{d_x^2}|z-x|^2\\
\le&2^{\alpha-1}d\frac{\|u_2\|_{L^\infty(B)}}{d_x^{1+\alpha}}|z-x|^{1+\alpha}.
\end{align*}

Concerning Condition \ref{cond_2}, consider $\partial\Omega_{2}=\Gamma\cup\partial\Omega$. If $y\in\Gamma$, Theorem \ref{22} yields the existence of a linear polynomial $P_y(z)$ such that
\begin{equation*}
\|P_y\|_{L^\infty(\Omega_{2})}+\|DP_y\|_{L^\infty(\Omega_{2})}\le E
\end{equation*}
and
\begin{equation*}
|u_2(z)-P_y(z)|\le E|z-y|^{1+\alpha},
\end{equation*}
for every $z\in\overline{\Omega_2}$, with $E\le C_0\|\psi\|_{C^{1,\alpha}(B_1')}\|g\|_{C^{0,\alpha}(\Gamma)}$, and $C_0=C_0(d,\alpha)>0$. If $y\in\partial\Omega\in C^\infty$, then, by classical boundary regularity for harmonic functions, $u_2\in C^{1,\alpha}\left(\overline{B\cap\Omega}\right)$, with $B:=B_r(y)$, for some $r>0$ sufficiently small. By Theorem \ref{23}, there exists a linear polynomial $P_y(z)$ such that
\begin{equation*}
|u_2(z)-P_y(z)|\le C_0|z-y|^{1+\alpha},
\end{equation*}
for every $z\in\overline{\Omega_2}$, for some $C_0=C_0(d,\alpha)>0$. This fact ensures that Condition \ref{cond_2} is in force.

Finally, we address Condition \ref{cond_3}. Fix $y\in\partial\Omega_{2}$, and let $P_y$ be the corresponding linear polynomial from Condition \ref{cond_2}. Clearly, the function $v=u_2-P_y$ is harmonic in $\Omega_{2}$, so it satisfies the interior estimates from Condition \ref{cond_1}.

Therefore, by Proposition \ref{24}, we have $u_2\in C^{1,\alpha}\left(\overline{\Omega_2}\right)$, and there exists a constant $C>0$, depending only on $d,\alpha,\Gamma$, such that 
\[
    \|u_2\|_{C^{1,\alpha}(\overline{\Omega_2})}\le C\|g\|_{C^{0,\alpha}(\Gamma)}.
\]

We close this section by discussing an alternative proof to Theorem \ref{21}; it appeared in \cite{DongP} and consists in re-casting the transmission problem in \eqref{1} as a Dirichlet problem driven by a PDE with a piecewise $C^{0,\alpha}$-regular right-hand side. Indeed, consider the problem of finding $w\in W^{1,2}(\Omega)$ satisfying
\begin{equation}\label{gatwick}
    \begin{cases}
        \Delta w=c&\hspace{.2in}\mbox{in}\hspace{.2in}\Omega_1\\
        w_\nu=g&\hspace{.2in}\mbox{on}\hspace{.2in}\partial\Omega_1,
    \end{cases}
\end{equation}
with 
\[
        \int_{\Omega_1}w{\rm d}x=0,
\]
where
\[
    c:=\frac{1}{\mathcal{H}^{d-1}(\Gamma)}\int_\Gamma g\,{\rm d}\mathcal{H}^{d-1}.
\]

The existence of solutions to \eqref{gatwick} follows from Sobolev inequalities and the Lax-Milgram Theorem. Also, usual results in elliptic regularity theory yield the existence of $\alpha\in (0,1)$ and $C>0$ such that $w\in C^{1,\alpha}(\Omega_1)$ with
\[
    \left\|w\right\|_{C^{1,\alpha}(\Omega_1)}\leq C\left\|g\right\|_{C^\alpha(\Gamma)}.
\]

Once the existence and regularity for the auxiliary function $w$ are available, we notice the solution $u$ to \eqref{1} satisfies
\begin{equation}\label{brickln}
    \begin{cases}
        \Delta u = -{\rm div}\left(\chi_{\Omega_1}Dw\right)+\chi_{\Omega_1}c&\hspace{.2in}\mbox{in}\hspace{.2in}\Omega\\
        u=0&\hspace{.2in}\mbox{on}\hspace{.2in}\partial\Omega.
    \end{cases}
\end{equation}

Once again, the Lax-Milgram Theorem ensures the existence of a weak solution to \eqref{brickln}. Concerning the regularity of $u$, an application of \cite[Corollary 2 and Remark 3]{Dong2012} yields the result. We notice the approach in \cite{DongP} accommodates more general elliptic operators in divergence form, as long as uniform ellipticity is available. Furthermore, it covers nonhomogeneous equations.

\section{A quasilinear fixed transmission problem}\label{sec_cpr}

In this section, we discuss an example of a quasilinear fixed transmission problem modelled after the $p$-Laplace operator. The gist of the section is to extend to a class of degenerate, variational models, some of the results in \cite{CSCS2021}, namely the regularity across the fixed transmission interface. The following material is based on the findings reported in \cite{Bianca-Pimentel-Urbano2023}.

We consider the following problem. Let $\Omega\subset\mathbb{R}^d$ be a bounded domain and fix $\Omega_1\Subset\Omega$. Define $\Omega_2:=\Omega\setminus\Omega_1$ and $\Gamma:=\partial\Omega_1$. Suppose that the interface $\Gamma$ is a $(d-1)$-manifold of class $C^1$. For a function $u:\overline{\Omega}\to\mathbb{R}$, we set
\begin{equation*}
u_1:=u\big|_{\overline{\Omega_1}}\hspace{.3in} \mbox{and}\hspace{.3in} u_2:=u\big|_{\overline{\Omega_2}}.
\end{equation*}
The normal derivatives of $u_1$ and $u_2$ on the interface are defined as
\[
    \frac{\partial u_i}{\partial\nu} := Du_i \cdot \nu, \quad i=1,2,
\]
where $\nu$ stands for the unit normal vector to $\Gamma$ inwards $\Omega_1$. 

The variational formulation of the problem relies on a functional driven by a nonlinear function $g:\mathbb{R}_0^+\to\mathbb{R}_0^+$, whose properties will be stated later. The PDE counterpart of the model is the quasilinear degenerate transmission problem consisting of finding a function $u:\overline{\Omega}\to\mathbb{R}$ such that
\begin{equation}\label{eq_stima118}
\begin{cases}
		\textnormal{div}\left(\frac{g\left(| D  u_1|\right)}{| D  u_1|} D  u_1\right)=0&\hspace{.2in}\mbox{in}\hspace{.2in}\Omega_1\\
		\vspace*{-0.3cm}\\
		\textnormal{div}\left(\frac{g\left(| D  u_2|\right)}{| D  u_2|} D  u_2\right)=0&\hspace{.2in}\mbox{in}\hspace{.2in}\Omega_2,\\
\end{cases}
\end{equation}
under the following boundary and interface conditions
\begin{equation}\label{eq_stima119}
\begin{cases}
u=0&\hspace{.2in}\mbox{on}\hspace{.2in}\partial\Omega\\
\vspace*{-0.3cm}\\
\frac{g(| D  u_1|)}{| D  u_1|}\frac{\partial u_1}{\partial\nu}-\frac{g(| D  u_2|)}{| D  u_2|}\frac{\partial u_2}{\partial\nu}=f&\hspace{.2in}\mbox{on}\hspace{.2in}\Gamma,
\end{cases}
\end{equation}
where the nonlinear function $f$ is given. The integrability properties of the latter will lead to two different regularity results for the solutions.

Regarding the nonlinear function $g\in C^1\left(\mathbb{R}^+_0\right)$, we suppose that 
	\begin{equation}
		g_0\le\frac{tg'(t)}{g(t)}\le g_1,\quad\forall t>0,
	\label{business class}
	\end{equation}
	for fixed constants $1\le g_0\le g_1$. In addition, we impose the monotonicity condition
	\begin{equation}\label{monicavitti}
		\bigg(\frac{g(|\xi|)}{|\xi|}\xi-\frac{g(|\zeta|)}{|\zeta|}\zeta\bigg)\cdot(\xi-\zeta)\ge C|\xi-\zeta|^p, \quad\forall\xi,\zeta\in\mathbb{R}^d,
	\end{equation}
for $p>2$, fixed though arbitrary, and $C>0$.

The work-horse of the theory is the case $g(t)=t^{p-1}$, with $p>2$, which turns \eqref{eq_stima118} into degenerate $p$-Laplace equations. A distinct example of nonlinearity within the scope of \eqref{business class}-\eqref{monicavitti} is
	\[
		g(t):=t^{p-1}\ln\left(a+t\right)^\alpha,
	\]
	for $p>2$, $a>1$ and $\alpha>0$. 
	
An important ingredient in the analysis is the primitive of $g$, 
	\begin{equation*}
		G(t)=\int_0^tg(s)\,{\rm d}s, \quad t \geq 0.
	\end{equation*}
We notice that $G$ is a Young function. That is, it is left-continuous and convex; we refer the reader to \cite[Definition 3.2.1]{at1}. Being a Young function, $G$ allows us to define a functional space. This is the subject of the next definition.

\begin{Definition}[Orlicz-Sobolev space]\label{notto}
Let $G$ be a Young function. The Orlicz-Sobolev space $W^{1,G}(\Omega)$ is the set of weakly differentiable functions $u\in W^{1,1}(\Omega)$ such that 
\[
	\int_\Omega G\left(|u(x)|\right){\rm d}x+\int_\Omega G\left(|Du(x)|\right){\rm d}x<\infty.
\]
\end{Definition}
Reasoning as in the usual setting, it is also possible to define $W^{1,G}_0(\Omega)$. Now we introduce a notion of solution to \eqref{eq_stima118}-\eqref{eq_stima119}.

\begin{Definition}\label{def_weaksol}
We say that $u\in W_0^{1,G}(\Omega)$ is a weak solution of \eqref{eq_stima118}-\eqref{eq_stima119} if 
\begin{equation}
\int_{\Omega}\frac{g\left(| D  u|\right)}{| D  u|} D  u\cdot D  v\,{\rm d}x=\int_{\Gamma}f v\,{\rm d}\mathcal{H}^{d-1},\hspace{.2in}\forall\,v\in W^{1,G}_0(\Omega).
\label{portia}
\end{equation}
\end{Definition}

\begin{Remark}
An important consideration on Definition \ref{def_weaksol} concerns the fact that the integrals in \eqref{portia} are well-defined. This is actually the case; let $u,v\in W^{1,G}(\Omega)$ and notice that
\begin{equation*}
	\int_{\Omega}\frac{g(|Du|)}{|Du|}Du\cdot Dv\,{\rm d}x<\infty.
\end{equation*}
Also $tg(t)\le CG(t)$, for $t\ge0$, because $g$ is increasing. Finally, $G(t+s)\le C\big(G(t)+G(s)\big)$, for $t,s\ge0$. As a result,
\begin{align*}
\bigg|\int_{\Omega}\frac{g(|Du|)}{|Du|}Du\cdot Dv\,{\rm d}x\bigg|\le&\int_{\Omega}g(|Du|)|Dv|\,{\rm d}x\\
\le&\int_{\Omega}g(|Du|+|Dv|)(|Du|+|Dv|)\,{\rm d}x\\
\le&C\int_{\Omega}G(|Du|+|Dv|)\,{\rm d}x\\
\le&C\int_{\Omega}G(|Du|)\,{\rm d}x+C\int_{\Omega}G(|Dv|)\,{\rm d}x\\
<&\infty.
\end{align*}
\end{Remark}

\begin{Remark}[On the connection between $W^{1,G}_0$ and $W^{1,p}_0$]
Let $u\in W_0^{1,G}(\Omega)$, and suppose that \eqref{monicavitti} holds true. Then we see that $u\in W_0^{1,p}(\Omega)$. Indeed, that inequality yields
\begin{align*}
\int_\Omega|Du|^p\,{\rm d}x\le&C\int_\Omega g(|Du|)|Du|\,{\rm d}x\\
\le&C\int_\Omega G(|Du|)\,{\rm d}x.
\end{align*}
\end{Remark}

\begin{Remark}[Existence and uniqueness of solutions]
The existence of a unique weak solution to \eqref{eq_stima118}-\eqref{eq_stima119} follows from approximation and monotonicity methods (see \cite{baroni2015}; see also \cite{Lieberman_1991}). An important characterization of the weak solutions to \eqref{eq_stima118} concerns the functional $I: W_0^{1,G}(\Omega)\to\mathbb{R}$ given by
\begin{equation}\label{stima117}
 I(u)=\int_{\Omega}G\left(| D u|\right)\,{\rm d}x-\int_{\Gamma}fu\,{\rm d}\mathcal{H}^{d-1}.
\end{equation}
Indeed, a weak solution to \eqref{eq_stima118}-\eqref{eq_stima119} is a global minimizer for $I$, whose first compactly supported variation yields \eqref{portia}.
\end{Remark}

Since our interest lies in the H\"older-types of moduli of continuity, we proceed by introducing functional spaces suitable to our analysis. These are the Campanato and Morrey spaces, which provide us with useful characterizations for the regularity estimates on \eqref{eq_stima118}-\eqref{eq_stima119}
\begin{Definition}[Campanato spaces]
We denote by $L_C^{p,\lambda}(\Omega;\mathbb{R}^d)$, with $1\le p<\infty$ and $\lambda\ge0$, the space of functions $u\in L^p(\Omega;\mathbb{R}^d)$ such that
\begin{equation*}
[u]_{L_C^{p,\lambda}(\Omega;\mathbb{R}^d)}^p=\sup_{x^0\in\Omega,\rho>0}\frac{1}{\rho^{\lambda}}\int_{\Omega\cap B(x^0,\rho)}|u-(u)_{\Omega\cap B(x^0,\rho)}|^p\,{\rm d}x<\infty.
\end{equation*}
\end{Definition}

\begin{Definition}[Morrey spaces]
We denote by $L_M^{p,\lambda}(\Omega;\mathbb{R}^d)$, with $1\le p<\infty$ and $\lambda\ge0$, the space of functions $u\in L^p(\Omega;\mathbb{R}^d)$ such that
\begin{equation*}
\|u\|_{L_M^{p,\lambda}(\Omega)}^p=\sup_{x^0\in\Omega,\rho>0}\frac{1}{\rho^{\lambda}}\int_{\Omega\cap B(x^0,\rho)}|u|^p\,{\rm d}x<\infty.
\end{equation*}
\end{Definition}

We close this section with two technical results. The first one relates a decay regime with a H\"older-type of inequality. 

\begin{Lemma}\label{lem_stima119}
Fix $R_0>0$ and let $\phi:[0,R_0]\to[0,\infty)$ be a non-decreasing function. Suppose there exist constants $C_1,\alpha,\beta>0$, and $C_2,\mu\ge0$, with $\beta<\alpha$, satisfying
\begin{equation*}
\phi(r)\le C_1\Big[\Big(\frac{r}{R}\Big)^{\alpha}+\mu\Big]\phi(R)+C_2R^{\beta},
\end{equation*}
for every $0<r\le R\le R_0$.
Then, for every $\sigma\le\beta$, there exists $\mu_0=\mu_0(C_1,\alpha,\beta,\sigma)$ such that, if $\mu<\mu_0$, for every $0<r\le R\le R_0$, we have
\begin{equation*}
\phi(r)\le C_3\Big(\frac{r}{R}\Big)^{\sigma}\big(\phi(R)+C_2R^{\sigma}\big),
\end{equation*}
where $C_3=C_3(C_1,\alpha,\beta,\sigma)>0$. Moreover,
\begin{equation*}
\phi(r)\le C_4r^{\sigma},
\end{equation*}
where $C_4=C_4(C_2,C_3,R_0,\phi(R_0),\sigma)$.
\end{Lemma}
Our second technical lemma relates solutions to the transmission problem with solutions to the homogeneous equation. Indeed, it provides a lower bound for the difference of their $G$-integrals in terms of $p$-norms.
\begin{Lemma}\label{stima146}
		Let $w\in W^{1,G}(B_R)$. Suppose \eqref{business class}--\eqref{monicavitti} is in force. Suppose further $h\in W_w^{1,G}(B_R)$ is a weak solution to
		\begin{equation*}
			\textnormal{div}\bigg(\frac{g(| D  h|)}{| D  h|} D  h\bigg)=0\quad\textnormal{in }B_R.
		\end{equation*}
Then there exists $C>0$ such that
		\begin{equation}\label{antonioni}
			\int_{B_{R}}G(| D  w|)-G(| D  h|)\,{\rm d}x\ge C\int_{B_{R}}| D (w-h)|^p\,{\rm d}x.
		\end{equation}
	\end{Lemma}
 For a proof of Lemma \ref{stima146} we refer the reader to \cite{Bianca-Pimentel-Urbano2023}

\subsection{Regularity estimates in ${\rm BMO}-$spaces}\label{subsec_bmo}

In the sequel, we consider $f\in L^\infty(\Omega)$. Our goal is to prove ${\rm BMO}-$regularity estimates for the gradient of solutions \emph{across the fixed transmission interface}. The following proposition is instrumental in our analysis.

\begin{Proposition}\label{stima136}
Let $h\in W^{1,G}(B_R)$ be a weak solution of
	\begin{equation*}
		\textnormal{div}\bigg(\frac{g(| D  h|)}{| D  h|} D  h\bigg)=0\quad\textnormal{in }B_R.
	\end{equation*}
Suppose \eqref{business class}-\eqref{monicavitti} hold true. Then there exist $C>0$ and $\alpha\in(0,1)$ such that, for every $r\in(0,R]$, we obtain
	\begin{equation*}
		\int_{B_r}| D  h-( D  h)_r|\,{\rm d}x\le C\Big(\frac{r}{R}\Big)^{d+\alpha}\int_{B_R}| D  h-( D  h)_R|\,{\rm d}x.
	\end{equation*}
\end{Proposition}
For a proof of Proposition \ref{stima136}, we refer the reader to \cite{baroni2015}.
\begin{Proposition}\label{stima161}
		Let $w\in W^{1,G}(B_R)$, and suppose $h\in W^{1,G}(B_R)$ is a weak solution of
		\begin{equation*}
			\textnormal{div}\bigg(\frac{g(| D  h|)}{| D  h|} D  h\bigg)=0\quad\textnormal{in }B_R.
		\end{equation*}
Suppose \eqref{business class}-\eqref{monicavitti} hold true. Then there exists $C>0$ such that, for every $0<r\le R$, one has
		\begin{align}
			\int_{B_r}| D  w-( D  w)_r|\,{\rm d}x\le& C\Big(\frac{r}{R}\Big)^{d+\alpha}\int_{B_R}| D  w-( D  w)_R|\,{\rm d}x \notag\\
			&+C\int_{B_R}| D  w- D  h|\,{\rm d}x \notag
		\end{align}
		for $\alpha\in(0,1)$ as in Proposition \ref{stima136}.
	\end{Proposition}
The proof of Proposition \ref{stima161} is technical in nature, and we choose to omit it here (see \cite[Proposition 2]{Bianca-Pimentel-Urbano2023} for details). The main result in \cite{Bianca-Pimentel-Urbano2023} reads as follows.

\begin{Theorem}[Gradient regularity in ${\rm BMO}-$spaces]\label{thm_ll}
Let $u\in W^{1,G}_0(\Omega)$ be a weak solution for the transmission problem \eqref{eq_stima118}--\eqref{eq_stima119}. Suppose \eqref{business class}--\eqref{monicavitti} are in force. Then $ Du\in {\rm BMO}_{\rm loc}(\Omega)$. In addition, for every $\Omega'\Subset\Omega$, 
\[
	\left\|Du\right\|_{{\rm BMO}(\Omega')}\leq C,
\]
where $C=C(d,\|f\|_{L^\infty(\Gamma)},{\rm diam}(\Omega),{\rm dist}(\Omega',\partial\Omega))>0$.
\end{Theorem}

The proof relies on an $L^1$-distance for the gradient of the solution $u$ and the $g$-harmonic function $h$ agreeing with $u$ in the Sobolev sense on the boundary of $B_R$. In fact, we aim at producing
\begin{equation}\label{eq_usb}
    \int_{B_R}\left|D\left(u-h\right)\right|{\rm d}x\leq CR^d.
\end{equation}
The former inequality follows from an involved combination of the minimality of $u$, properties of $g$-harmonic functions, Lemma \ref{stima146} and standard results, such as the Trace Theorem and the Poincar\'e Inequality. 

Then \eqref{eq_usb} builds upon Proposition \ref{stima161} to produce
				\begin{equation*}
					\int_{B_r}| D  u-( D  u)_r|\,{\rm d}x\le C\Big(\frac{r}{R}\Big)^{d+\alpha}\int_{B_R}| D  u-( D  u)_R|\,{\rm d}x+CR^d,
				\end{equation*}
				for every $0<r\le R$. 
    Finally, an application of Lemma \ref{lem_stima119}, implies
\begin{equation*}
    \int_{B_r}| D  u-( D  u)_r|\,{\rm d}x\le Cr^d, \quad\forall r\in(0,R],
\end{equation*}
which completes the proof.

At this point one resorts to embedding results for borderline spaces to obtain a modulus of continuity for the solution $u$ in $C^{0,{\rm Log-Lip}}-$spaces. See \cite[Theorem 3]{Cianchi1996}. In fact, we have the following corollary.

\begin{Corollary}[Log-Lipschitz continuity estimates]\label{cor_ll}
Let $u\in  W_0^{1,G}(\Omega)$ be a weak solution for \eqref{eq_stima118}-\eqref{eq_stima119}. Suppose \eqref{business class}--\eqref{monicavitti} hold true. Then $u\in C^{0,{\rm Log-Lip}}_{\rm loc}(\Omega)$. Also, for every $\Omega'\Subset\Omega$, 
\[
	\left\|u\right\|_{C^{0,{\rm Log-Lip}}(\Omega')}\leq C\left(\left\|u\right\|_{L^\infty(\Omega)}+\left\|f\right\|_{L^\infty(\Gamma)}\right),
\]
where $C=C(p,d,{\rm diam}(\Omega),{\rm dist}(\Omega',\partial\Omega))>0$.
\end{Corollary}

To complete this section we mention the case of unbounded interface data. Suppose $f\in W^{1,p'+\varepsilon}(\Omega)$, where $\varepsilon>0$ depends on $p$ and the dimension $d$. Hence, it is possible to prove the following regularity result in, H\"older spaces, for the solutions to \eqref{eq_stima118}-\eqref{eq_stima119}.

\begin{Theorem}\label{thm_c0alpha}
Let $u$ be a weak solution to the interface problem \eqref{eq_stima118}-\eqref{eq_stima119}, and suppose \eqref{business class}-\eqref{monicavitti} are in force. Let $2<p<d$, and $\varepsilon>0$ be such that
\[
	\frac{d-p}{p-1}<\varepsilon<d-\frac{p}{p-1}.
\]
Suppose further that $f\in W^{1,p'+\varepsilon}(\Omega)$. Then $u\in C^{0,\alpha}_{\textnormal{loc}}(\Omega)$, with
\[
	\alpha=1-\frac{d}{p+\varepsilon(p-1)},
\]
and estimates are available.
\end{Theorem}

\begin{Remark}[Local boundedness]
We note that arguing along the same ideas as those presented in \cite{Serrin_1964}, one can prove that solutions to \eqref{eq_stima118}-\eqref{eq_stima119} are locally bounded (at the level of the functions and their gradients). Although we omit this result and its proof in the present manuscript, the reader may find the details in \cite[Section 3]{Bianca-Pimentel-Urbano2023}.
\end{Remark}

\begin{Remark}[Potential estimates and the $p$-Laplace operator]
For the case $g(t):=t^{p-1}$, one can use \cite[Corollary 1, item (C9)]{KM2014c} to derive the conclusion of Theorem \ref{thm_ll}. Because the interface $\Gamma$ is locally of class $C^1$, it follows from the inequality
\[
	\int_{B_r\cap \Gamma}f{\rm d}\mathcal{H}^{d-1}\leq Cr^{d-1};
\]
see, for instance, \cite[Proposition 3.5]{Maggi2012}.
\end{Remark}

\section{Free transmission problems}\label{subsec_ussr}

A variant of fixed transmission problems concerns models whose interfaces are solution-dependent. In this scenario, the models appear in the context of free boundary problems. As a consequence, an additional structure arises as part of the unknown; namely, the free interface. 

We present this class of problems in the context of fully nonlinear elliptic equations, which are intrinsically non-variational and rely on the notion of viscosity solutions. We start by introducing the problem.

\subsection{A fully nonlinear free transmission problem: the setting}\label{subsec_fnlsetting}

Let $\Omega\subset\mathbb{R}^d$ be a bounded domain and consider $u\in C(\Omega)$. We define $\Omega^\pm(u)$ as 
\[
    \Omega^+(u):=\left\lbrace x\in\Omega\,|\,u(x)>0 \right\rbrace
\]
and
\[
    \Omega^-(u):=\left\lbrace x\in\Omega\,|\,u(x)<0 \right\rbrace.
\]

Now, let $0<\lambda\leq \Lambda$ be fixed constants. We let $F_1,F_2:S(d)\to\mathbb{R}$ be $(\lambda,\Lambda)$-elliptic operators. That is,
\[
    \lambda\|N\|\leq F_i(M+N)-F_i(M)\leq \Lambda\|N\|,
\]
for $i\in\left\lbrace1,2\right\rbrace$, and every symmetric matrix $M,N\in S(d)$, with $N\geq 0$.
We are interested in the problem 
\begin{equation}\label{eq_ashmolean}
    \begin{cases}
        F_1(D^2u)=f&\quad \mbox{in}\quad \Omega^+(u)\\
        F_2(D^2u)=f&\quad \mbox{in}\quad \Omega^-(u).
    \end{cases}
\end{equation}
The (free) interface $\Gamma_u$ in this case is the topological boundary of the set where $u\neq 0$. That is,
\[
    \Gamma_u:=\partial\left\lbrace u>0\right\rbrace\cup \partial\left\lbrace u<0\right\rbrace.
\]
We consider two classes of solutions to \eqref{eq_ashmolean}: viscosity solutions (both in the continuous and the $L^p$-senses) and $L^p$-strong solutions. For the specifics on those notions, we refer the reader to \cite{CIL,CCKS}; see also \cite[Chapter 17]{Gilbarg-Trudinger}. It is worth noticing that when considering $L^p$-strong solutions one derives information on the transmission condition from the Sobolev properties of the solutions.

Because an $L^p$-strong solution is a function in $W^{2,d}(B_1)$, we have $Du=0$ almost everywhere in $\left\lbrace u=0\right\rbrace$. Hence, it is natural to equip \eqref{eq_ashmolean} with the condition
\[
    \left|Du\right|=0\hspace{.3in}\mbox{on}\hspace{.3in}\partial\left\lbrace u>0\right\rbrace\cup \partial\left\lbrace u<0\right\rbrace.
\]

A number of difficulties arise in the analysis of \eqref{eq_ashmolean}. First, the existence of solutions is a non-trivial matter. Indeed, because the diffusion process is discontinuous with respect to the solutions, the dependence of the equation on $u$ is unknown. In particular, there is no a priori reason to expect it to be monotone. Hence, the equation governing the free transmission problem may lack properness and, therefore, a comparison principle. 

In addition, the regularity of the solutions can no longer rely on the geometry of the transmission interface since such an object is unknown. As a consequence, the approach developed in \cite{SoriaCarro-Stinga} does not yield information on \eqref{eq_ashmolean}. In the sequel, we present recent developments bypassing those difficulties and extending the existence and regularity program to \eqref{eq_ashmolean}.

\subsection{The existence of solutions}\label{subsec_tobeornottobe}

We consider the Dirichlet problem 
\begin{equation}\label{eq_towerroom}
    \begin{cases}
        F_1(D^2u)\chi_{\{u>0\}}+F_2(D^2u)\chi_{\{u<0\}}=f&\hspace{.05in}\mbox{in}\hspace{.05in}\left(\Omega^+(u)\cup\Omega^-(u)\right)\cap\Omega\\
        u=g&\hspace{0.05in}\mbox{on}\hspace{0.05in}\partial\Omega,
    \end{cases}
\end{equation}
where $\Omega\subset\mathbb{R}^d$ is a bounded domain satisfying a uniform exterior one condition, $f\in L^p(\Omega)$, for $p>p_0$, and $g\in C(\partial\Omega)$ is a given boundary condition. We recall that $d/2\leq p_0=p_0(\Lambda/\lambda,d)$ is the exponent such that $(\lambda,\Lambda)$-elliptic equations with right-hand side in $L^p$, with $p>p_0$, are entitled to the Aleksandrov-Bakelman-Pucci maximum principle. 

\begin{Remark}[Escauriaza exponent $d/2\leq p_0=p_0(\Lambda/\lambda,d)$]\label{rem_escorial}
    The connection of the exponent $p_0$ with fully nonlinear elliptic equations of the form $F(D^2u,x)=f$ appears in the work of Luis Escauriaza \cite{escauriaza}. In that paper, the author proves a Harnack inequality of the form
    \begin{equation}\label{eq_baldwin}
        \sup_{x\in B_{r/2}}u(x)\leq C\left(\inf_{x\in B_{r/2}}u(x)+r^{2-d/q}\left\|f\right\|_{L^p(B_r)}\right),
    \end{equation}
    provided $u\geq 0$ solves $F(D^2u,x)=f$ and $f\in L^p(B_1)$, for $p>p_0$. This inequality follows from the improved integrability of the Green's function associated with the linearization of $F$. Such improved integrability of the Green's function was established in the work of Eugene Fabes and Daniel Stroock \cite{fabes_stroock1984}. In the derivation of \eqref{eq_baldwin}, one notices that Escauriaza's exponent $p_0$ is the conjugate of the improved integrability available for the Green's function of a $(\lambda,\Lambda)$-elliptic operator. 
\end{Remark}

We present the following theorem.

\begin{Theorem}[Existence of viscosity solutions]\label{thm_existence1}
Let $\Omega\subset\mathbb{R}^d$ be a bounded domain satisfying a uniform exterior cone condition. Let $F_i:S(d)\to\mathbb{R}$ be $(\lambda,\Lambda)$-elliptic operators, for $i\in\{1,2\}$, and some fixed constants $0<\lambda\leq \Lambda$. Suppose $g\in C(\partial\Omega)$ and $f\in L^p(\Omega)$, for some $p>p_0$.
Then there exists an $L^p$-viscosity solution $u\in C(\overline\Omega)$ to \eqref{eq_towerroom}.
\end{Theorem}

Before detailing the proof of Theorem \ref{thm_existence1}, we mention that further conditions on the operators $F_1$ and $F_2$ allow us to obtain qualitative information on the solution to \eqref{eq_towerroom} whose existence follows from the theorem. The first condition concerns the local proximity of the operators $F_1$ and $F_2$; we suppose there exist constants $K,\tau>0$ such that 
\begin{equation}\label{eq_assump1}
    \left|F_1(M)-F_2(M)\right|\leq K+\tau\left\|M\right\|,
\end{equation}
for every $M\in S(d)$. Condition \eqref{eq_assump1} unlocks a $C^{1,\alpha}$-estimate for the $L^p$-viscosity solution in Theorem \ref{thm_existence1}. 

The second condition we impose on $F_i$ yields the existence of $L^p$-strong solutions to \eqref{eq_towerroom}. This is a convexity-type of assumption and reads as follows. Suppose there exists $L,\sigma>0$, and a \emph{convex} $(\lambda,\Lambda)$-elliptic operator $F:S(d)\to\mathbb{R}$ such that 
\begin{equation}\label{eq_expressabs}
    \left|F_i(M)-F(M)\right|\leq L+\sigma\left\|M\right\|,
\end{equation}
for every $M\in S(d)$.

\begin{Theorem}[Existence of $C^{1,\alpha}$-regular solutions]\label{thm_existence2}
Let $\Omega\subset\mathbb{R}^d$ be a bounded domain satisfying a uniform exterior cone condition. Let $F_i:S(d)\to\mathbb{R}$ be $(\lambda,\Lambda)$-elliptic operators, for $i\in\{1,2\}$, and some fixed constants $0<\lambda\leq \Lambda$. Suppose $g\in C(\partial\Omega)$ and $f\in L^p(\Omega)$, for some $p>p_0$. Suppose further that \eqref{eq_assump1} also holds and $p>d$; let $\alpha\in(0,1)$ satisfy
 \[
 	\alpha<\alpha_0\hspace{.2in}\mbox{and}\hspace{.2in}\alpha\leq 1-\frac{d}{p},
 \]
 where $\alpha_0\in(0,1)$ corresponds to the $C^{1,\alpha_0}$-regularity available for the solutions to $G=0$ for any $(\lambda,\Lambda)$-elliptic operator $G$. Then there exists $\beta_0=\beta_0(d,p,\lambda,\Lambda,\alpha)>0$ such that, if the parameter $\tau>0$ in \eqref{eq_assump1} satisfies $\tau\leq \beta_0$, then $u\in C^{1,\alpha}_{\rm{loc}}(\Omega)$ and, for every 
 $\Omega'\Subset\Omega$, we have
 \begin{equation*}\label{eq_c1aest}
\|u\|_{C^{1,\alpha}(\Omega')}\leq C\left(1+|F_1(0)|+|F_2(0)|+\|f\|_{L^p(\Omega)}
+\|g\|_{L^\infty(\partial\Omega)}\right),
\end{equation*}
where $C=C(\alpha,d,p,\lambda,\Lambda,K,\tau,{\rm diam}(\Omega),{\rm dist}(\Omega',\partial\Omega))$.
\end{Theorem}

We notice Theorem \ref{thm_existence2} not only ensures the existence of an $L^p$-viscosity solution to the free transmission problem but also provides a regularity estimate for this object. Of course, this is not a regularity result in the sense it applies only to the solution whose existence stems from the theorem. 

If we require $F_1$ and $F_2$ to be close, in the sense of \eqref{eq_assump1}, but also to satisfy a proximity regime with respect to a convex $(\lambda,\Lambda)$-operator $F$ as in \eqref{eq_expressabs}, the conclusion of Theorem \ref{thm_existence2} improves substantially. In fact, instead of a $C^{1,\alpha}$-regular $L^p$-viscosity solution, it is possible to establish the existence of a solution in $W^{2,p}(\Omega)$. As a consequence, under \eqref{eq_expressabs}, it is possible to establish the existence of $L^p$-strong solutions for the free transmission problem. This is the content of the next theorem.

\begin{Theorem}[Existence of $L^p$-strong solutions]\label{thm_existence3}
Let $\Omega\subset\mathbb{R}^d$ be a bounded domain satisfying a uniform exterior cone condition. Let $F_i:S(d)\to\mathbb{R}$ be $(\lambda,\Lambda)$-elliptic operators, for $i\in\{1,2\}$, and some fixed constants $0<\lambda\leq \Lambda$. Suppose $g\in C(\partial\Omega)$ and $f\in L^p(\Omega)$, for some $p>p_0$. There exists $\beta_0=\beta_0(d,p,\lambda,\Lambda,\alpha)>0$ such that, if the parameter $\sigma>0$ in \eqref{eq_expressabs} satisfies $\sigma\leq \beta_0$, then \eqref{eq_towerroom} has an $L^p$-strong solution $u\in W^{2,p}(\Omega)\cap C(\overline{\Omega)}$. In  addition, for every 
 $\Omega'\Subset\Omega$, 
 \begin{equation*}\label{eq_c2aest}
\|u\|_{W^{2,p}(\Omega')}\leq C\left(1+|F_1(0)|+|F_2(0)|+\|f\|_{L^p(\Omega)}
+\|g\|_{L^\infty(\partial\Omega)}\right),
\end{equation*}
where $C=C(\alpha,d,p,\lambda,\Lambda,L,\sigma,{\rm diam}(\Omega),{\rm dist}(\Omega',\partial\Omega))$.
\end{Theorem}

The main difficulty in establishing the existence of viscosity solutions for \eqref{eq_towerroom} stems from the dependence of the operator on the zeroth order term $u$. In fact, the lack of properness rules out standard formulations of the comparison principle. As a result, Perron's method is no longer available. 

However, a two-parameters regularization of the PDE in \eqref{eq_towerroom} turns the free transmission problem into an equation holding in the entire domain $\Omega$, with no explicit dependence on the solution $u$. At this level, the comparison principle is available, and one can prove the existence of global barriers. Perron's method yields the existence of viscosity solutions. Then one proceeds by applying a fixed-point argument at the level of the functional parameter. Secondly, one sends the remaining parameter to zero; stability of viscosity solutions recovers a solution to \eqref{eq_towerroom} and the proof is complete. 

In what follows, we introduce some of the ingredients in those arguments. The main idea concerns the regularization of the equation in \eqref{eq_towerroom}. Let $v\in C(\overline{\Omega})$ agree with the Dirichlet data $g$ on the boundary. Fix $\varepsilon>0$, arbitrary. Define a function $g_\varepsilon^v:\Omega\to\mathbb{R}$ as
\[
    g_\varepsilon^v(x):=\max\left(\min\left(\frac{v(x)+\varepsilon}{2\varepsilon},1\right),0\right)\hspace{.2in}\mbox{in}\hspace{.2in}\Omega,
\]
with $g_\varepsilon^v\equiv 0$ in $\mathbb{R}^d\setminus\Omega$. Notice that $g_\varepsilon^v(x)=1$ in $\left\{v>\varepsilon\right\}$, whereas $g_\varepsilon^v(x)=0$ in $\left\{v<-\varepsilon\right\}$. We consider the convolution of $g_\varepsilon^v$ with a standard mollifying kernel $\eta_\varepsilon$ to obtain
\[
    h_\varepsilon^v:=g_\varepsilon^v\ast\eta_{\varepsilon}.
\]
Once $h_\varepsilon^v$ is available, we introduce the auxiliary operator $G_\varepsilon^v:\Omega\times S(d)\to\mathbb{R}$ given by
\[
    G_\varepsilon^v(x,M):=h_\varepsilon^v(x)F_1(M)+(1-h_\varepsilon^v(x))F_2(M).
\]

A straightforward computation implies three properties of the operator $G_\varepsilon^v$. First, because $F_1$ and $F_2$ are $(\lambda,\Lambda)$-elliptic, we conclude $G_\varepsilon^v$ is also a $(\lambda,\Lambda)$-elliptic operator. One can also prove the existence of a constant $K_\varepsilon^v>0$ such that 
\[
    \left|G_\varepsilon^v(x,M)-G_\varepsilon^v(y,M)\right|\leq K_\varepsilon^v|x-y|\left(1+\left\|M\right\|\right),
\]
for every $x,y\in\Omega$ and every $M\in S(d)$. That is, we conclude that $G_\varepsilon^v$ is Lipschitz-continuous in the space-variable, locally uniformly in $M$.

A fundamental information on $G_\varepsilon^v$ concerns its connection with the conditions in \eqref{eq_assump1} and \eqref{eq_expressabs}. Indeed, the definition of $G_\varepsilon^v$ preserves those conditions, up to an adjustment in the constants. If \eqref{eq_assump1} holds, we have
\[
    \left|G_\varepsilon^v(x,M)-G_\varepsilon^v(y,M)\right|\leq 2\left(K+\tau\left\|M\right\|\right),
\]
for every $x,y\in\Omega$ and every $M\in S(d)$. Also, if \eqref{eq_expressabs} is in force, we get
\[
    \left|G_\varepsilon^v(x,M)-F(M)\right|\leq L+\sigma\left\|M\right\|,
\]
for every $x\in\Omega$ and $M\in S(d)$. For a proof of those properties, we refer the reader to \cite[Lemma 1]{Pimentel-Swiech2022}.

Once $G_\varepsilon^v$ is well defined, we consider the Dirichlet problem
\begin{equation}\label{eq_fatwa}
\left\{\begin{array}{rcll}
G_\varepsilon^v(x,D^2u_\varepsilon^v)&=&f\quad\mbox{in}\quad\Omega\\
        u_\varepsilon^v&=&g\quad\mbox{on}\quad\partial\Omega.
\end{array}\right.
\end{equation}
In \cite{Pimentel-Swiech2022}, the authors observe the comparison principle is available for the equation in \eqref{eq_fatwa}. In addition, they establish the existence of global sub and supersolutions $\underline u$ and $\overline u$ and notice that such functions are independent of $v$ and $\varepsilon>0$. Perron's method then implies the existence of an $L^p$-viscosity solution $u_\varepsilon^v$ to \eqref{eq_fatwa}. Moreover, $\underline u\leq u_\varepsilon^v\leq \overline u$.

At this point, one considers two objects. First, a subset $B\subset C(\overline\Omega)$ defined as
\[
    B:=\left\lbrace v\in C(\overline\Omega)\,|\;\underline u\leq v\leq \overline u\right\rbrace.
\]
Then we introduce an operator $T$ defined on $B$ as follows: for $v\in B$, one solves \eqref{eq_fatwa} to obtain the unique solution $u_\varepsilon^v$. One sets $Tv:=u_\varepsilon^v$. These ingredients satisfy certain properties. In fact, $B$ is a close and convex subset of $C(\overline\Omega)$. Also, $T$ maps $B$ to itself and $T(B)$ is a precompact set in $C(\overline\Omega)$. Finally, one notices that $T:B\to B$ is continuous. 

As a consequence, an application of the Schauder Fixed Point Theorem yields the existence of an $L^p$-viscosity solution to 

\begin{equation}\label{eq_fatwa2}
\left\{\begin{array}{rcll}
G_\varepsilon^u(x,D^2u_\varepsilon)&=&f\quad\mbox{in}\quad\Omega\\
u_\varepsilon&=&g\quad\mbox{on}\quad \partial\Omega.
\end{array}\right.
\end{equation}

Finally, to prove the existence of an $L^p$-viscosity solution to \eqref{eq_towerroom}, one takes the limit $\varepsilon\to 0$ and resorts to the stability of viscosity solutions. These steps lead to the conclusion of Theorem \ref{thm_existence1}. To establish Theorem \ref{thm_existence2} under \eqref{eq_assump1}, one can argue through well-understood arguments, such as in \cite[Theorem 8.3]{Caffarelli_Cabre1995} or \cite[Theorem 2.1]{Swiech1997}. Concerning Theorem \ref{thm_existence3}, one notices that \eqref{eq_expressabs} frames the problem in the context of classical $W^{2,p}$-regularity theory for fully nonlinear elliptic equations; see \cite[Chapter 7]{Caffarelli_Cabre1995}.

Once we have discussed the existence of the solutions to \eqref{eq_towerroom}, we proceed by addressing their regularity estimates. This is the subject of the next section.

\subsection{Regularity for fully nonlinear free transmission problems}\label{subsec_morrison}

Once the existence of $L^p$-viscosity and $L^p$-strong solutions has been addressed in the literature, the natural question concerns their regularity. This is the subject of \cite{Pimentel-Santos2023}, whose general lines we discuss in what follows.

In that work, the authors examine $L^p$-strong solutions to \eqref{eq_ashmolean} and establish two classes of regularity estimates. Under a near-convexity condition as in \eqref{eq_expressabs}, they show that strong solutions as locally of class $C^{1,{\rm Log-Lip}}$. Furthermore, under an additional (pointwise) condition on the density of the negative phase, they prove that solutions satisfy a (pointwise) quadratic growth regime.

Once again, the main difficulty in the analysis of \eqref{eq_ashmolean} is the discontinuity of the operator with respect to the solution. Indeed, standard results such as the Harnack inequality or the maximum principle are not available. To bypass the lack of usual ideas and methods the argument in \cite{Pimentel-Santos2023} relates \eqref{eq_ashmolean} with two viscosity inequalities, \emph{holding in the entire domain}, for which the theory is available.

In fact, let $u\in W^{2,d}(\Omega)$ be an $L^d$-strong solution to \eqref{eq_ashmolean}. Hence, it satisfies
\[
    \min\left\lbrace F_1(D^2u(x)),F_2(D^2u(x))\right\rbrace\leq 1 \hspace{.3in}\mbox{a.e. in}\hspace{.1in}\Omega
\]
and 
\[
    \max\left\lbrace F_1(D^2u(x)),F_2(D^2u(x))\right\rbrace\geq -1 \hspace{.3in}\mbox{a.e. in}\hspace{.1in}\Omega.
\]

One then explores the connection between $L^d$-strong and $L^d$-viscosity solutions, as detailed in \cite{CCKS}. It follows from \cite[Lemma 2.5]{CCKS} that $u$ is an $L^d$-viscosity solution to
\begin{equation}\label{eq_paris}
    \begin{cases}
        \min\left\lbrace F_1(D^2u),F_2(D^2u)\right\rbrace\leq 1&\hspace{.3in}\mbox{in}\quad\Omega\\
        \max\left\lbrace F_1(D^2u),F_2(D^2u)\right\rbrace\geq -1&\hspace{.3in}\mbox{in}\quad\Omega.
    \end{cases}
\end{equation}
We notice, however, the ingredients in \eqref{eq_paris} are continuous. Hence $u$ is indeed a $C$-viscosity solution to this pair of inequalities, and one immediately concludes that $u\in C^\alpha_{\rm loc}(\Omega)$, with estimates. This fact yields compactness for any family of solutions to \eqref{eq_paris}. Also, the operators governing these inequalities are entitled to stability results for viscosity solutions, in the spirit of \cite[Proposition 4.11]{Caffarelli_Cabre1995}.

\begin{Theorem}[Local $C^{1,\llip}$-regularity]\label{thm_main1}
Let $u\in W^{2,d}_{\rm{loc}}(B_1)$ be a strong solution to \eqref{eq_ashmolean}. Suppose $F_1$ and $F_2$ are convex, $(\lambda,\Lambda)$-elliptic operators. Suppose further that \eqref{eq_expressabs} is in force. There exists $0<\beta_0=\beta_0(d,\lambda,\Lambda)\ll1$ such that, if $\sigma\leq\beta_0$, we have $u\in C^{1,\llip}_{\rm{loc}}(B_1)$, and there exists $C>0$ such that
\[
	\sup_{x\in B_r(x_0)}\left|u(x)-u(x_0)-Du(x_0)\cdot(x-x_0)\right|\leq Cr^2\ln\frac{1}{r},
\]
for every $x_0\in\Omega$ and $r>0$ satisfying $B_r(x_0)\Subset\Omega$. In addition, $C=C(d,\lambda,\Lambda,\|u\|_{L^\infty(\Omega)})$.
\end{Theorem}

The proof of Theorem \ref{thm_main1} follows from approximation methods, along with the ideas introduced in the work of Luis Caffarelli \cite{Caffarelli1989}; see also \cite{Caffarelli1988_Milano}. For a detailed argument, we refer the reader to \cite{Pimentel-Santos2023}.

An interesting question concerning the solutions of \eqref{eq_ashmolean} arises from Theorem \ref{thm_main1}. Indeed, one may enquire into the further conditions on the problem that would switch the regularity regime from the borderline scenario $C^{1,{\rm Log-Lip}}$ to that of $C^{1,1}$-regularity estimates. 

A condition that allows us to improve the information in Theorem \ref{thm_main1} regards the density of the negative phase in a vicinity of a given point $x_0\in\Omega$ on the free boundary. We proceed by introducing some ingredients.

The free boundary associated with \eqref{eq_ashmolean} is the set 
\[
    \Gamma(u):=\left(\partial\Omega^+(u)\cup\partial\Omega^-(u)\right)\cap \Omega.
\]
A point $x_0\in \Gamma(u)$ may fall within three distinct categories. First, it can be a \emph{one-phase point}, in the sense that it lies on the free boundary between the positive or the negative phase and the set $\{u=0\}$. More rigorously, we say that $x_0\in\Gamma(u)$ is a one-phase point if
\[
	x_0\in\left(\partial\Omega^\pm(u)\setminus\partial\Omega^\mp(u)\right)\cap \Omega.
\]
When examining a one-phase free boundary point, the regularity analysis amounts to the study of an obstacle problem, which has been well-understood and documented; see, for instance, \cite{Lee1998}. 

Another situation refers to points $x_0\in\Gamma(u)$ sitting on a portion of the free boundary separating the positive and the negative phases. These are called \emph{two-phase points} and satisfy
\[
	x_0\in\left(\partial\Omega^+(u)\cap\partial\Omega^-(u)\right)\cap \Omega.
\]
The study of the regularity of the solutions or the free boundary properties around two-phase points is completely open in the context of \eqref{eq_ashmolean}. 

However, among such points, we highlight $x^*\in\Gamma(u)$ for which 
\begin{equation}\label{eq_branchp}
	|B_r(x^*)\cap\{u=0\}|>0,
\end{equation}
for every $0<r\ll1$. A two-phase point satisfying \eqref{eq_branchp} is called a \emph{branch point} for \eqref{eq_ashmolean}. 
We denote with $\Gamma_{\rm{BR}}(u)\subset \Gamma(u)$ the set of branch points for the solution $u$.

A further ingredient in the analysis of quadratic growth for the solutions to \eqref{eq_ashmolean} is the density of the negative phase. Given a point $x_0\in \Gamma(u)$, we define $V_r(x_0,u)$ as
\begin{equation}\label{eq_volume0}
	V_r(x_0,u) := \frac{{\rm Vol}\left(B_r(x_0)\cap \Omega^-(u)\right)}{r^d};
\end{equation}
the quantity $V_r(x_0,u)$ amounts to the density of the negative phase around a free boundary point and plays a critical role in the regularity theory for the solutions. Indeed, in \cite{Pimentel-Santos2023} the authors suppose $F_1$ and $F_2$ to be convex $(\lambda,\Lambda)$-elliptic operators, positively homogeneous of degree one. In addition, if there exist a constant $C_0>0$ and $x_0\in \Gamma(u)$ such that 
\[
    V_r(x_0,u)\leq C_0,
\]
they prove
\[
	\sup_{x\in B_r(x_0)}\left|u(x)\right|\leq Cr^2,
\]
for every $0<r\ll1$, provided $C_0$ can be taken \emph{arbitrarily small.}

Were the smallness condition available locally along $\Gamma(u)$, the quadratic growth would imply $C^{1,1}$-regularity estimates for the solutions to \eqref{eq_ashmolean}. However, such a condition would imply the negative phase is empty in such a neighbourhood, trivializing the problem. Instead of aiming at a local $C^{1,1}$-regularity result, the authors consider branch points $x^*\in\Gamma_{\rm BR}(u)$ and prove the following theorem.

\begin{Theorem}[Quadratic growth at branch points]\label{thm_main2}
Let $u\in W_{\rm loc}^{2,d}(B_1)$ be a strong solution to \eqref{eq_ashmolean}. Suppose the operators $F_1$ and $F_2$ are convex, positively homogeneous of degree one, and $(\lambda,\Lambda)$-elliptic. Let $x^*\in \Gamma_{\rm{BR}}(u)$ be such that $V_r(x_0,u)\leq C_0$. Then there exists $0<\overline{C_0}\ll1$ such that, if $C_0\leq \overline{C_0}$, one can find a universal constant $C>0$ for which
\[
	\sup_{x\in B_r(x^*)}\left|u(x)\right|\leq Cr^2,
\]
for every $0<r\ll1$.
\end{Theorem}

We present the general lines of the proof of Theorem \ref{thm_main2} in the sequel. It relies on a dyadic analysis combined with the maximum principle and a scaling argument. The latter uses the $L^\infty$-norms of $u$ as a normalization factor, introducing an additional dependence on this quantity into the estimates. 

The strategy used in \cite{Pimentel-Santos2023} is inspired by arguments first launched in \cite{Caffarelli-Karp-Shahgholian2000} in the analysis of a free boundary arising in the Pompeiu problem. In that paper, the diffusion process is driven by the Laplace operator. The fully nonlinear counterpart of \cite{Caffarelli-Karp-Shahgholian2000} appeared in \cite{Lee-Shahgholian2001}. Here, the authors identify that convexity and homogeneity of degree one are the precise conditions allowing to switch from the linear to the nonlinear setting. 

To detail the proof of Theorem \ref{thm_main2}, we start by defining a subset of the natural numbers related to upper bounds for $u$ in dyadic balls. Fix $x^*\in \Gamma_{\rm{BR}}(u)\cap \Omega$. The maximal subset of $\mathbb{N}$ whose elements $j$ are such that
\begin{equation}\label{eq_16ju}
	\sup_{x\in B_{2^{-j-1}}(x^*)}|u(x)| \geq \frac{1}{16}\sup_{x\in B_{2^{-j}}(x^*)}|u(x)|
\end{equation}
is denoted with $\mathcal{M}(x^*,u)$. 

\begin{Proposition}\label{prop_quadratic}
Let $u\in W_{\rm loc}^{2,d}(B_1)$ be a strong solution to \eqref{eq_ashmolean}. Suppose the operators $F_1$ and $F_2$ are convex, positively homogeneous of degree one, and $(\lambda,\Lambda)$-elliptic. Let $x^*\in \Gamma_{\rm{BR}}(u)$ be such that
\begin{equation}\label{eq_smalldensity}
	V_{2^{-j}}(x^*,u) < C_0,
\end{equation}
for every $j\in \mathcal{M}(x^*,u)$, where $C_0\leq \overline{C_0}$. Then
\[
	\sup_{x\in B_{2^{-j}}(x^*)} |u(x)| \leq \frac{1}{C_0}2^{-2j}, \hspace{.4in} \forall j \in \mathcal{M}(x^*,u).
\]
\end{Proposition}

For a proof of Proposition \ref{prop_quadratic}, we refer the reader to \cite[Proposition 3]{Pimentel-Santos2023}. Once the quadratic growth holds in the set $\mathcal{M}(x^*,u)$, one extends this fact to the natural numbers $\mathbb{N}$. This is the subject of the next proposition.

\begin{Proposition}\label{prop_quadnaturals}
Let $u\in W_{\rm loc}^{2,d}(B_1)$ be a strong solution to \eqref{eq_ashmolean}. Suppose the operators $F_1$ and $F_2$ are convex, positively homogeneous of degree one, and $(\lambda,\Lambda)$-elliptic. Let $x^*\in \Gamma_{\rm{BR}}(u)$ be such that
\begin{equation}\label{eq_smalldensity2}
	V_{2^{-j}}(x^*,u) < C_0,
\end{equation}
for every $j\in \mathcal{M}(x^*,u)$, where $C_0\leq \overline{C_0}$. Then
\[
	\sup_{x\in B_{2^{-j}}(x^*)} |u(x)| \leq \frac{4}{C_0}2^{-2j}, \hspace{.4in} \forall j \in \mathbb{N}.
\]
\end{Proposition}

For the detail in the proof of Proposition \ref{prop_quadnaturals}, we refer the reader to \cite[Proposition 4]{Pimentel-Santos2023}. Finally, a discrete-to-continuous argument builds upon Proposition \ref{prop_quadnaturals} to complete the proof of Theorem \ref{thm_main2}.

Given $0<r\ll1$, let $j\in\mathbb{N}$ be such that $2^{-(j+1)}\leq r< 2^{-j}$. Hence,
\[
	\begin{split}
		\sup_{B_r}|u(x)| \leq \sup_{B_{2^{-j}}}|u(x)| \leq C\left[\left(\frac{1}{2}\right)^{j+1-1}\right]^2 \leq Cr^2,
	\end{split}
\]	
and the proof is complete.

\bigskip

{\small \noindent{\bf Acknowledgments.} All authors are partially supported by the Centre for Mathematics of the University of Coimbra (UIDB/00324/2020, funded by the Portuguese Government through FCT/MCTES). EP is partially supported by FAPERJ (grants E26/200.002/2018 and E26/201.390/2021). JMU is partially supported by the King Abdullah University of Science and Technology (KAUST).}

\end{document}